\let\confidentialstring=\relax
\def\confidentialstring{%
\textit{IEEE Trans. Automat. Control}, vol. 56, 2011, accepted,
version of \today.%
}
\let\headnote=\relax
\def\mykeywords{Discrete abstraction, symbolic model, nonlinear
  system, symbolic control, motion planning, formal verification,
  polyhedral over-approximation, attainability, attainable
  set}

\documentclass[journal,a4paper,onecolumn,nofonttune,12pt,draftclsnofoot]{IEEEtran}
\edef\baselinestretch{1}
\def\confidentialstring{%
This work has been submitted to the IEEE for possible publication.
Copyright may be transferred without notice, after which this
version may no longer be accessible.%
}
\def\headnote{%
Please refer to
\href{http://www.reiszig.de/gunther/pubs/i11abs.html}{author's homepage}
for the latest version of this work and for a link to the definite
publication, respectively.
You may also find a {Bib\TeX} entry at this website.%
}
\def\mykeywords{Discrete abstraction, symbolic model, nonlinear
  system, symbolic control, motion planning, formal verification,
  polyhedral over-approximation, attainability, attainable
  set\ifCLASSOPTIONonecolumn; MSC: Primary, 93C10; Secondary, 93C55,
  93C57, 93C15, 93B03\fi}
\def\confidentialstring{%
This work has been accepted for publication in the
\emph{IEEE Trans. Automatic Control}.
Copyright may be transferred without notice, after which this
version may no longer be accessible.%
}
\def\headnote{%
Please refer to
\href{http://www.reiszig.de/gunther/pubs/i11abs.html}{author's homepage}
and to
\href{http://ieeexplore.ieee.org/}{IEEE Xplore}
for the definite publication.
You may also find a {Bib\TeX} entry at the former website.%
}
\def\myname{Gunther Rei\ss ig}

\def\mytitle{Computing abstractions of nonlinear systems}

\usepackage{cite}
\usepackage{psfrag}
\usepackage{algorithmic}
\usepackage{graphicx}
\usepackage{xcolor}

\usepackage[cmex10]{amsmath}
\let\ORGforeignlanguage\foreignlanguage
\def\foreignlanguage#1{\lowercase{\ORGforeignlanguage{#1}}}
\makeatletter
\def\MakeUppercase#1{#1}%
\def\markboth#1#2{\def\leftmark{\@IEEEcompsoconly{\sffamily}\MakeUppercase{#1}}%
\def\rightmark{\@IEEEcompsoconly{\sffamily}\MakeUppercase{#2}}}
\makeatother
\ifCLASSOPTIONonecolumn% not for final layout

\fi
\usepackage{gunther2e}

\ifx
\hypersetup\undefined\def\href#1#2{\texttt{#2}}
\else
\hypersetup{
pdftitle={\mytitle},
pdfauthor={Gunther Reissig, http://www.reiszig.de/gunther/},%
pdfsubject={\confidentialstring{} \headnote},
pdfkeywords={\mykeywords}
}
\fi

\begin{document}%
\makeatletter
%%%%%%%%%%%%%%%%%%%% correct proof environment %%%%%%%%%%%%%%%%%%%%%%%
% had been overwritten by IEEEtran, to the effect that QED commands
% did not work properly anymore
\renewenvironment{proof}[1][\proofname]{\par
  \pushQED{\qed}%
  \normalfont \topsep6\p@\@plus6\p@\relax
  \trivlist
  \item[\hskip\labelsep
        \itshape
    #1\@addpunct{.}]\ignorespaces
}{%
  \popQED\endtrivlist\@endpefalse
}
%%%%%%%%%%%%%%%%%%%% define headline %%%%%%%%%%%%%%%%%%%%%%%%%%%%%%%%%
\markboth{\myname\hspace*{\fill}\mytitle\hspace*{\fill}}%
{\myname\hspace*{\fill}\mytitle\hspace*{\fill}}%
%%%%%%%%%%%%%%%%%%%%%%%%%%%%%%%%%%%%%%%%%%%%%%%%%%%%%%%%%%%%%%%%%%%%%%
\makeatother

\title{%
%\raisebox{1cm}[0pt][0pt]{\makebox[0pt][l]{\normalsize\confidentialstring}}%
\mytitle}

\author{\myname%
\ifCLASSOPTIONonecolumn$^\ast$\else,~\IEEEmembership{Senior Member, IEEE}\fi%
\thanks{%
\ifCLASSOPTIONonecolumn$^\ast$\fi%
Universit{\"a}t Kassel,
Fachbereich 16 - Elektrotechnik/Informatik,
Regelungs- und Systemtheorie,
Wilhelmsh{\"o}her Allee 73,
D-34121 Kassel,
Germany,
\url{http://www.reiszig.de/gunther/}%
}% <-this % stops a space
\ifCLASSOPTIONonecolumn%
\thanks{\confidentialstring{} \headnote}%
\fi%
}

\maketitle

\begin{abstract}
Sufficiently accurate finite state models, also called symbolic models
or discrete abstractions, allow one to apply fully automated methods,
originally developed for purely discrete systems, to formally reason
about  continuous and hybrid systems, and to design finite state
controllers that provably enforce predefined specifications.
We present a novel algorithm to compute such finite state models
for nonlinear discrete-time and sampled systems which depends on
quantizing the state space using polyhedral cells, embedding these
cells into suitable supersets whose attainable sets are convex, and
over-approximating attainable sets by intersections of supporting
half-spaces. We prove a novel recursive description of
these half-spaces and propose an iterative procedure to compute them
efficiently. We also provide new sufficient conditions for the convexity
of attainable sets which imply the existence of the aforementioned
embeddings of quantizer cells. Our method yields highly accurate
abstractions and applies to nonlinear systems under mild assumptions,
which reduce to sufficient smoothness in the case of sampled
systems. Its practicability in the design of discrete
controllers for nonlinear continuous plants under state and control
constraints is demonstrated by an example.
\end{abstract}


\section{Introduction}
\label{s:intro}

In recent years, there has been a
growing interest in using finite state models for the analysis and
synthesis of continuous and hybrid systems
\cite{%
Hsu87,%
AlurHenzingerLafferrierePappas00,%
KoutsoukosAntsaklisStiverLemmon00,%
DellnitzJunge02,%
DingLiZhou02,%
Schroeder03,%
BlankeKinnaertLunzeStaroswiecki06,%
Osipenko07,%
Tabuada09%
}.
This interest has been stimulated by safety critical applications
\cite{TomlinPappasSastry98},
inherent limits of continuous feedback control
\cite[Sections 5.8-5.10]{Sontag98},
increasingly complex control objectives
% possible addition: FainekosGirardKressGazitPappas09
\cite{LaValle06},
and the necessity to cope with the effects of coarse quantization
\cite{MatveevSavkin09}. A sufficiently
accurate finite state model, also called a \begriff{symbolic model} or
\begriff{discrete abstraction}, would allow one to apply fully
automated methods, originally developed for purely discrete systems
\cite{KumarGarg95,ClarkeGrumbergPeled99,CassandrasLafortune08}, to
formally reason about the original system, and to design finite state
controllers that provably enforce predefined specifications
\cite{%
AlurHenzingerLafferrierePappas00,%
KoutsoukosAntsaklisStiverLemmon00,%
DellnitzJunge02,%
DingLiZhou02,%
Schroeder03,%
BlankeKinnaertLunzeStaroswiecki06,%
Osipenko07,%
Tabuada09%
}.
Obtaining such abstractions constitutes a
challenging problem, which has only been satisfactorily solved for
special cases.

Under the name \begriff{symbolic dynamics}, finite state models of
continuous systems had already been a well-established mathematical
tool \cite{LindMarcus95} when the concept
appeared in the engineering literature
\cite{CrandallChandiramaniCook66,Wang68}.
Much of the subsequent research
has been devoted to systems whose
continuous-valued dynamics is
\label{11958:ii}
linear.
Methods for nonlinear systems have
been systematically studied since around 1980; see
\cite{Hsu87,DellnitzJunge02,DingLiZhou02,Schroeder03,BlankeKinnaertLunzeStaroswiecki06,Osipenko07,Tabuada09}.
In the earliest such approach \cite{HsuGuttaluZhu82},
attainable sets are approximated
by means of trajectories emanating from a finite set of initial
points, hence the name \begriff{sampling method} \cite{HsuGuttaluZhu82}.
This method has been successfully applied to a
variety of problems
\cite{Hsu87,DellnitzJunge02,DingLiZhou02,Gruene02b,Schroeder03,BlankeKinnaertLunzeStaroswiecki06,Osipenko07,Hsu95,GrueneJunge07},
including symbolic control of sampled systems
\cite{Hsu85,GrueneJunge07}.
An extension
allows for rigorous over-approximation of attainable sets
\cite{Junge99},
and thus, for the  computation of abstractions.
Over the years, a large number of alternatives to the sampling method
have been proposed, which represent a variety of compromises between
approximation accuracy, practicability, rigor, and computational
complexity
\cite{Tabuada08,PolaTabuada09,
MoorRaisch02,MalerBatt08,BermanHalaszKumar07,GirardPappas05,
JaulinWalter97,StursbergKowalewskiEngell00,AlthoffStursbergBuss09,TazakiImura09,
Tiwari08b,
ChutinanKrogh03,MitchellBayenTomlin05,GrueneJunge07,GrueneMueller08,
Broucke98,CainesWei98,StiverKoutsoukosAntsaklis01}.

In the present paper, we aim at computing abstractions for nonlinear
discrete-time systems of the form
\begin{equation}
\label{e:TimeDiscreteAutonomousControlSystem:G}
x_{k+1} = G(x_k,u_k),
\end{equation}
where the state $x$ takes values in a subset of $\mathbb{R}^n$, and $u$ is an input
signal which is assumed to take its values in some finite set $U$.
If \ref{e:TimeDiscreteAutonomousControlSystem:G} arises from a
continuous-time system
\begin{equation}
\label{e:TimeContinuousAutonomousControlSystem:F}
\dot x = F(x,v)
\end{equation}
under sampling, its right hand side $G$ may not be
explicitly given. Our results will still apply as
we will formulate hypotheses to be verified and computations to be
performed directly in terms of the right hand side $F$ of
\ref{e:TimeContinuousAutonomousControlSystem:F}.

The approach we follow involves quantizing the state space of
\ref{e:TimeDiscreteAutonomousControlSystem:G}
with the help of a finite covering $C$ of $\mathbb{R}^n$ whose elements
we call \begriff{cells}
\cite{%
Hsu87,AlurHenzingerLafferrierePappas00,KoutsoukosAntsaklisStiverLemmon00,DellnitzJunge02,%
DingLiZhou02,Schroeder03,BlankeKinnaertLunzeStaroswiecki06,%
Osipenko07,Tabuada09%
}.
The system \ref{e:TimeDiscreteAutonomousControlSystem:G}
is supplemented with a \begriff{quantizer} $Q$ which assigns to any
state $x$ of \ref{e:TimeDiscreteAutonomousControlSystem:G} the
collection of those cells in $C$ that contain $x$,
$Q(x) = \Menge{\Delta \in C}{x \in \Delta}$. That is, a pair
% $(u,\Delta) \colon \mathbb{Z}_+ \to U \times C$
$(u,\Delta)$
of an input signal $u_0, u_1, \dots$ and an output signal
$\Delta_0, \Delta_1, \dots$ could possibly be generated by the
\begriff{quantized system} composed of
\ref{e:TimeDiscreteAutonomousControlSystem:G} and the
non-deterministic output relation
\begin{equation}
\label{e:TimeDiscreteAutonomousControlSystem:Q}
\Delta_k \in Q(x_k)
\end{equation}
iff there exists a sequence $x_0, x_1, \dots$ such that
\ref{e:TimeDiscreteAutonomousControlSystem:G} and
\ref{e:TimeDiscreteAutonomousControlSystem:Q} hold for all
non-negative integers $k$%
\footnote{The symbols $u$, $x$ and $\Delta$ are used to
denote elements of $U$, $\mathbb{R}^n$ and $C$, respectively, as well
as signals taking their values in these sets.}.
The collection of such pairs $(u,\Delta)$ is
called the \begriff{behavior} of the quantized system
\ref{e:TimeDiscreteAutonomousControlSystem:G},\ref{e:TimeDiscreteAutonomousControlSystem:Q}
\cite{Willems89}.

The input alphabet $U$ and the output alphabet $C$ of the quantized system
\ref{e:TimeDiscreteAutonomousControlSystem:G},\ref{e:TimeDiscreteAutonomousControlSystem:Q}
are both finite. Control problems for
\ref{e:TimeDiscreteAutonomousControlSystem:G},\ref{e:TimeDiscreteAutonomousControlSystem:Q}
can still be challenging to solve, especially if the system
\ref{e:TimeDiscreteAutonomousControlSystem:G} is nonlinear and the
specification involves constraints or is otherwise complex. In
contrast, controllers (or \begriff{supervisors}) for finite automata
are generally straightforward to design
\cite{KumarGarg95,CassandrasLafortune08}, which raises the question of
whether controllers for the quantized system
\ref{e:TimeDiscreteAutonomousControlSystem:G},\ref{e:TimeDiscreteAutonomousControlSystem:Q}
can be obtained by solving auxiliary control problems for automata that
approximate the behavior of
\ref{e:TimeDiscreteAutonomousControlSystem:G},\ref{e:TimeDiscreteAutonomousControlSystem:Q}.
As it turns out, this strategy is feasible if the approximation
is both conservative and sufficiently precise,
e.g.
\cite{KoutsoukosAntsaklisStiverLemmon00,Moor99,MoorRaisch99,MoorDavorenAnderson02}. That
is, the automaton must be capable of generating any signal in the behavior of
\ref{e:TimeDiscreteAutonomousControlSystem:G},\ref{e:TimeDiscreteAutonomousControlSystem:Q},
and the set of spurious signals should be small.
In other words, the said strategy requires a \begriff{discrete abstraction},
by which we mean a superset of the behavior of
\ref{e:TimeDiscreteAutonomousControlSystem:G},\ref{e:TimeDiscreteAutonomousControlSystem:Q}
that can be realized by a finite
automaton, and this abstraction should be as accurate as possible.

\begin{figure}[!t]
\begin{minipage}{\ifCLASSOPTIONonecolumn.49\else.99\fi\linewidth}
\centering
\psfrag{x0}[][]{\textcolor{red}{$x_t$}}%$x_0$
\psfrag{x1}[r][r]{\textcolor{red}{$x_{t+1}$}}%$x_1$}}
\psfrag{x2}[][]{\textcolor{red}{$x_{t+2}$}}%$x_2$}}
\psfrag{xl}[r][r]{\textcolor{red}{$x_{t+N}$}}%$x_k$}}
\psfrag{G(.,u0)}[][]{\small$G(\cdot,u_t)$}
\psfrag{G(.,u1)}[][]{\small$G(\cdot,u_{t+1})$}
\psfrag{G(.,u2)}[][]{\small$G(\cdot,u_{t+2})$}
%\psfrag{F}[][]{$\small\varphi(T,\cdot)$}
\psfrag{F}[][]{}
\psfrag{O0}[][]{$\Delta_{t}$}%$\Delta_{0}$}
\psfrag{O1}[][]{$\Delta_{t+1}$}%$\Delta_{1}$}
\psfrag{O2}[][]{$\Delta_{t+2}$}%$\Delta_{2}$}
\psfrag{Ol}[][]{$\Delta_{t+N}$}%$\Delta_{k}$}
\includegraphics[width=\ifCLASSOPTIONonecolumn3.5in\else.99\linewidth\fi,keepaspectratio]{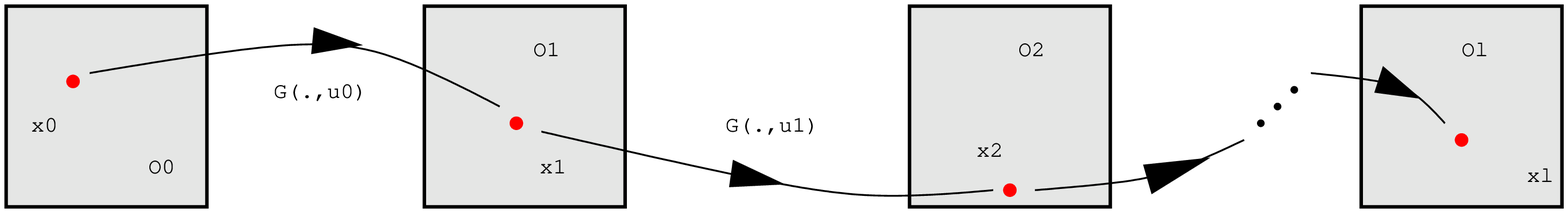}\\
(a)
\end{minipage}
\ifCLASSOPTIONonecolumn\hspace*{\fill}\else\\[.5ex]\fi
\begin{minipage}{\ifCLASSOPTIONonecolumn.48\else.98\fi\linewidth}
\centering
\psfrag{x0}[lb][lb]{\ifCLASSOPTIONonecolumn\small\fi\textcolor{red}{$x_0$}}
\psfrag{x1}[lt][lt]{}%\ifCLASSOPTIONonecolumn\scriptsize\fi\textcolor{red}{$x_1$}}
\psfrag{x2}[l][l]{}%\ifCLASSOPTIONonecolumn\small\fi\textcolor{red}{$x_2$}}
\psfrag{x1t}[lb][lb]{\ifCLASSOPTIONonecolumn\scriptsize\fi\textcolor{red}{$x_1$}}%\tilde x_1$}}
\psfrag{x2t}[l][l]{}%\ifCLASSOPTIONonecolumn\small\fi\textcolor{red}{$\tilde x_2$}}
\psfrag{D0}[bl][bl]{\ifCLASSOPTIONonecolumn\footnotesize\else\small\fi$\Delta_0$}
\psfrag{D1}[bl][bl]{\ifCLASSOPTIONonecolumn\footnotesize\else\small\fi$\Delta_1$}
\psfrag{D2}[tl][tl]{}%\ifCLASSOPTIONonecolumn\footnotesize\else\small\fi$\Delta_2$}
\psfrag{D2t}[bl][bl]{\ifCLASSOPTIONonecolumn\footnotesize\else\small\fi$\Delta_2$}%$\widetilde\Delta_2$}
\psfrag{Gb}[b][b]{\ifCLASSOPTIONonecolumn\tiny\else\scriptsize\fi$G(\cdot,u_0)$}
\psfrag{Gt}[t][t]{\ifCLASSOPTIONonecolumn\tiny\else\scriptsize\fi$G(\cdot,u_1)$}
\psfrag{1}[t][t]{\ifCLASSOPTIONonecolumn\footnotesize\else\small\fi$1$}
\psfrag{2}[t][t]{\ifCLASSOPTIONonecolumn\footnotesize\else\small\fi$2$}
\psfrag{xm1}[t][t]{\ifCLASSOPTIONonecolumn\footnotesize\else\small\fi$-1$}
\psfrag{xm2}[t][t]{\ifCLASSOPTIONonecolumn\footnotesize\else\small\fi$-2$}
\psfrag{y1}[r][r]{\ifCLASSOPTIONonecolumn\footnotesize\else\small\fi$1$}
\psfrag{ym1}[r][r]{\ifCLASSOPTIONonecolumn\footnotesize\else\small\fi$-1$}
\includegraphics[width=.625\linewidth,keepaspectratio]{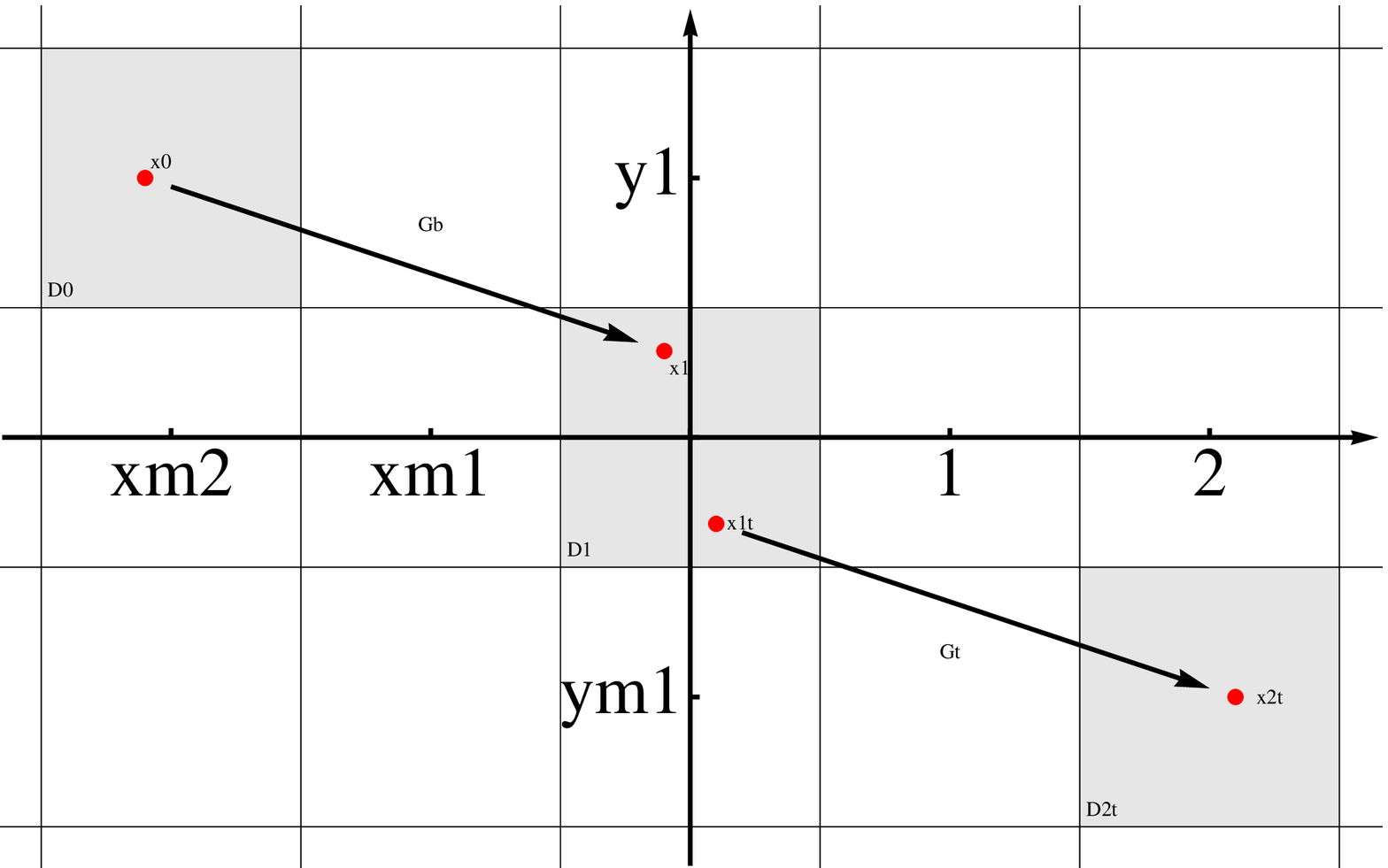}\\
(b)
\end{minipage}
\caption{\label{fig:CharNCompleteHull.GQ.pointwise}
(a) Illustration of
% condition \ref{e:Mcondition}, where $x_{\tau} = \psi(\tau,x_0,u)$.
consistency with the dynamics of the quantized system
\ref{e:TimeDiscreteAutonomousControlSystem:G},\ref{e:TimeDiscreteAutonomousControlSystem:Q}
on finite time intervals.
(b) Quantizer cells may serve as states of automata
realizations of abstractions of memory span $1$, where transitions are
defined by condition \ref{e:intro:consistent:N=1}.
}
\end{figure}

One way to prescribe the accuracy of an abstraction is to restrict the
extent by which its signals are allowed to violate the dynamics of
\ref{e:TimeDiscreteAutonomousControlSystem:G},\ref{e:TimeDiscreteAutonomousControlSystem:Q}.
While, by definition, signals in the behavior of the quantized system
\ref{e:TimeDiscreteAutonomousControlSystem:G},\ref{e:TimeDiscreteAutonomousControlSystem:Q}
are consistent with the dynamics of
\ref{e:TimeDiscreteAutonomousControlSystem:G},\ref{e:TimeDiscreteAutonomousControlSystem:Q}
at all times, a common class of abstractions require consistency only
on finite time intervals \cite{Moor99,MoorRaisch99}.
Such an abstraction contains
any pair $(u,\Delta)$ that fulfills the following condition, in which
the \begriff{memory span} $N \geq 1$ determines the accuracy
\cite{Willems89}: For any non-negative integer $t$ there are
states $x_t, \dots, x_{t+N}$ such that
\ref{e:TimeDiscreteAutonomousControlSystem:G} and
\ref{e:TimeDiscreteAutonomousControlSystem:Q} hold for all
$k \in \{t,t+1,\dots,t+N-1 \}$ and $k \in \{t, t+1, \dots, t+N \}$,
respectively. See \ref{fig:CharNCompleteHull.GQ.pointwise}(a).

In the case that $N = 1$, consistency of $(u,\Delta)$ is equivalent to
the existence of a sequence $x_0, x_1, \dots$ for which
\begin{equation}
\label{e:intro:consistent:N=1}
x_k \in \Delta_k\text{\ \ and\ \ }G(x_k, u_k) \in \Delta_{k+1}
\end{equation}
hold for all $k$.
Hence, the cells in the covering $C$ may serve as states of an
automaton realization of the abstraction, where the occurrence of an
input symbol $u_k \in U$ enables a transition from
$\Delta_k \in C$ to $\Delta_{k+1} \in C$ iff there is a state $x_k$ of
\ref{e:TimeDiscreteAutonomousControlSystem:G} such
that \ref{e:intro:consistent:N=1} holds. Obviously, that automaton
will be capable of generating any pair of signals $u$ and $\Delta$ in
the behavior of
\ref{e:TimeDiscreteAutonomousControlSystem:G},\ref{e:TimeDiscreteAutonomousControlSystem:Q}. The
fact that it will generally also generate spurious signals is
illustrated in \ref{fig:CharNCompleteHull.GQ.pointwise}(b).
If \ref{e:TimeDiscreteAutonomousControlSystem:G}
requires the sign of the second component of the state to be constant, then
the sequence $\Delta_0, \Delta_1, \Delta_2, \dots$ of cells generated
by the automaton is not consistent with the dynamics of
\ref{e:TimeDiscreteAutonomousControlSystem:G},\ref{e:TimeDiscreteAutonomousControlSystem:Q}.
In contrast, consistency of ($u,\Delta)$ for $N > 1$ requires, amongst
other conditions, that \ref{e:intro:consistent:N=1} holds with
$x_1 = G(x_0, u_0)$, which rules out the spurious signal
$\Delta_0, \Delta_1,\Delta_2, \dots$ of
\ref{fig:CharNCompleteHull.GQ.pointwise}(b).
Indeed, increasing the memory span $N$ generally results in more
accurate abstractions.

In this paper we shall present a novel algorithm to compute
abstractions of finite but otherwise arbitrary memory span that builds
on a well-known reformulation of consistency on finite time intervals
in terms of attainable sets \cite{Moor99,MoorRaisch99}, on a new
method to compute polyhedral over-approximations of the latter, and on
new results that guarantee the convexity of attainable sets of
\ref{e:TimeDiscreteAutonomousControlSystem:G} and
\ref{e:TimeContinuousAutonomousControlSystem:F}.

In our approach, polyhedral quantizer cells are embedded into suitable
supersets whose attainable sets under the dynamics of
\ref{e:TimeDiscreteAutonomousControlSystem:G} are convex for the
duration of $N$ time steps, where $N$ is the memory span of the
abstraction that is being computed.
That convexity requirement permits us to over-approximate attainable sets
by intersections of supporting half-spaces, and the latter are
obtained from systems of linear equations derived from
\ref{e:TimeDiscreteAutonomousControlSystem:G}. The number of
half-spaces needed can be quite large, especially
if the memory span exceeds $1$. We present a novel recursive
description of these half-spaces and propose an iterative procedure to
compute them efficiently.

The existence of the aforementioned embeddings of quantizer cells
is, in fact, the essential requirement for our method to apply. The
results in this paper not only allow verification of that requirement
when a particular quantizer is given, but they also show how to meet it
using sufficiently small but otherwise arbitrary polyhedral cells.
We use strongly convex supersets of quantizer cells, and the error by
which we over-approximate attainable sets depends quadratically on the
size of the cells. 
Application of our earlier results \cite{i07Convex,i07MMAR} on
ellipsoidal supersets would have led to linear error bounds.
Thus, the accuracy of the computed abstractions is improved if a
particular quantizer is given. Alternatively, fewer and larger cells
may be used, which reduces the computational effort to compute
abstractions and also reduces the complexity of controllers designed
on the basis of the latter.
These results are obtained under mild
assumptions on the right hand side $G$ of
\ref{e:TimeDiscreteAutonomousControlSystem:G}, which reduce to
sufficient smoothness in the case of sampled systems.

The remaining of this paper is organized as follows.
The next section introduces basic notation and terminology.
In Section \ref{s:computation} we present our algorithm for the
computation of abstractions, prove its correctness, and analyze its
computational complexity.
Section \ref{s:convex} is devoted to our results on the convexity of
attainable sets.
In Section \ref{s:ex}, practicability of our
approach in the design of discrete controllers for nonlinear
continuous plants under state and control constraints is demonstrated
by an example. We also present computational results on how
the computational effort of our approach grows with the problem size.

\section{Preliminaries}
\label{s:prelims}

\subsection{Basic notation}

$\mathbb{R}$ and $\mathbb{Z}$ denote the sets of real numbers and
integers, respectively, $\mathbb{R}_{+}$ and $\mathbb{Z}_{+}$, their
subsets of non-negative elements,
and $\mathbb{N} = \mathbb{Z}_{+} \setminus \{ 0 \}$.
$\intcc{a,b}$, $\intoo{a,b}$,
$\intco{a,b}$, and $\intoc{a,b}$
denote closed, open and half-open, respectively,
intervals with end points $a$ and $b$, e.g.
$\intco{0,\infty} = \mathbb{R}_{+}$.
$\intcc{a;b}$, $\intoo{a;b}$,
$\intco{a;b}$, and $\intoc{a;b}$ stand for discrete intervals, e.g.
$\intcc{a;b} = \intcc{a,b} \cap \mathbb{Z}$.

For any sets $A$ and $B$,
$f \colon A \to B$ denotes a map of $A$ into $B$,
and $B^A$ is the set of all such maps.
%$f|_C$ is the restriction of $f$ to $C$,
Operations involving subsets of $\mathbb{R}^n$ are
defined pointwise \cite[Appendix A]{Valentine64}, e.g.
% $\alpha M \defas \Menge{\alpha y}{y \in M}$,
$\Delta + \Delta' \defas \Menge{\omega + \omega'}{\omega \in \Delta, \omega' \in \Delta'}$
and
$\varphi(\intcc{0,t},\Delta)
\defas
\Menge{\varphi(\tau,\omega)}%
{\tau \in \intcc{0,t}, \omega \in \Delta}$ if
$\Delta,\Delta' \subseteq \mathbb{R}^n$,
$\varphi \colon \mathbb{R} \times \mathbb{R}^n \to \mathbb{R}^n$, and
$t \in \mathbb{R}$.
% $\alpha \in \mathbb{R}$ and $M,N \subseteq \mathbb{R}^n$.
% %

$C^k$ denotes the class of $k$ times continuously differentiable maps,
and $C^{k,1}$, the class of maps in $C^k$ with (locally) Lipschitz-continuous
$k$th derivative.

\subsection{Behaviors}
\label{ss:behave}

Given an arbitrary set $W$ called \begriff{signal alphabet}, any
subset $B \subseteq W^{\mathbb{Z}_+}$ is a
\begriff{behavior} on $W$
\cite{Willems89,Moor99,MoorRaisch99,MoorDavorenAnderson02,i09HSCC}.
Hence, elements of $B$ are infinite sequences
$w \colon \mathbb{Z}_{+} \to W$, which we call \begriff{signals}.
We denote the value of the signal $w$ at time $k$ by $w_k$.
The \begriff{backward $\tau$-shift} $\sigma^{\tau}$ is defined by
$
\left(\sigma^{\tau} w\right)_k = w_{\tau+k}
$.
The \begriff{restriction} of $B$ to $I \subseteq \mathbb{Z}_+$,
$B|_{I}$, is defined by $B|_{I} \defas \Menge{w|_{I}}{w \in B}$.
$B$ is \begriff{time-invariant} if $\sigma^1 B \subseteq B$.
% , where $\sigma = \sigma^1$.
$B$ is \begriff{$N$-complete}, or equivalently, $B$ has
\begriff{memory span} $N$, if $N \in \mathbb{Z}_+$ and
$
B
=
\Menge{w \in W^{\mathbb{Z}_+}}{
\forall_{\tau \in \mathbb{Z}_+}\;
(\sigma^{\tau} w)|_{\intcc{0;N}} \in B|_{\intcc{0;N}}
}
$.
A superset $B'$ of a behavior $B \subseteq W^{\mathbb{Z}_+}$ is called
an \begriff{abstraction} of $B$, and $B'$ is additionally called
\begriff{discrete} if it can be realized by a finite automaton.

\subsection{Discrete-time systems}
\label{ss:SystemClass}

In \ref{e:TimeDiscreteAutonomousControlSystem:G} with right hand side
$G \colon X \times U \to X$,
$u \colon \mathbb{Z}_{+} \to U$ represents an input signal and
$x \colon \mathbb{Z}_{+} \to X$, a state signal.
% The input alphabet $U$ will often be assumed to be finite and may be a
% singleton, in which case \ref{e:TimeDiscreteAutonomousControlSystem:G}
% is an autonomous ordinary difference equation.
A \begriff{trajectory} of
\ref{e:TimeDiscreteAutonomousControlSystem:G} is a sequence
$(x,u) \colon \mathbb{Z}_+ \to X \times U$
for which \ref{e:TimeDiscreteAutonomousControlSystem:G} holds for all
$k \in \mathbb{Z}_+$.
The collection of such trajectories, which is a
subset of $\left( X \times U \right)^{\mathbb{Z}_+}$, is
called the \begriff{behavior} of
\ref{e:TimeDiscreteAutonomousControlSystem:G}.
The \begriff{general solution}
$\psi \colon \mathbb{Z}_+ \times X \times U^{\mathbb{Z}_+} \to X$
of \ref{e:TimeDiscreteAutonomousControlSystem:G} is the map defined by
the requirement that $(\psi(\cdot,x_0,u),u)$ is a trajectory of
\ref{e:TimeDiscreteAutonomousControlSystem:G} and
$\psi(0,x_0,u) = x_0$.
% the maximal solution of
% \ref{e:TimeDiscreteAutonomousControlSystem:G} for which
% $\psi(0,x_0,u) = x_0$, for all $x_0 \in \mathbb{R}^n$ and all
% $u \colon \mathbb{Z}_+ \to U$.
Of course, it is not necessary to specify $u$
on the whole time axis, so we may write
\[
\psi(k,x_0,u_0,\dots,u_{k' - 1}) \defas \psi(k,x_0,u|_{\intco{0,k'}}) \defas \psi(k,x_0,u)
\]
whenever $k' \geq k$.
% The following hypothesis is introduced for later reference.
We will often assume the following.

\begin{hypothesis}
\label{h:GdiffeomorphismC1}
$X \subseteq \mathbb{R}^n$ is open and
$G \colon X \times U \to X$ is such that for all $u \in U$, the map
$G(\cdot,u)$ is a $C^1$-diffeomorphism onto an open subset
of $X$.
\end{hypothesis}

We will say \begriff{\ref{h:GdiffeomorphismC1} is fulfilled with smoothness
$C^K$} if \ref{h:GdiffeomorphismC1} holds with $G(\cdot,u)$ of class
$C^K$ rather than merely $C^1$.

\section{Computation of abstractions}
\label{s:computation}

In the present section, we shall develop an
efficient algorithm for computing abstractions
of \ref{e:TimeDiscreteAutonomousControlSystem:G} that
can be represented by finite automata.
Intuitively, our approach is that of successively expanding the behavior
of \ref{e:TimeDiscreteAutonomousControlSystem:G} and may be seen to
comprise four approximation steps:
State space quantization,
approximation by the smallest discrete abstraction,
approximation by a collection of convex programs,
approximation by a collection of linear programs.

The purpose of state space quantization is to conservatively approximate
\ref{e:TimeDiscreteAutonomousControlSystem:G} using a finite signal
alphabet, which is an important prerequisite for a finite automaton
approximation. Unfortunately, $1$-completeness of the behavior of
\ref{e:TimeDiscreteAutonomousControlSystem:G} is usually lost, and in
general the behavior of the quantized system is not $N$-complete for
any $N$.
In order to reintroduce $N$-completeness, which is sufficient for a
finite automaton representation to exist, we approximate the quantized
system again, this time by the smallest discrete abstraction of memory
span $N$. A problem with the latter abstraction is that it may only be
computed exactly for special cases of both systems
\ref{e:TimeDiscreteAutonomousControlSystem:G} and state space
quantizations. Two more approximation steps yield further abstractions,
which are both $N$-complete and characterized in terms of
computationally tractable problems. Specifically, we first replace
each quantizer cell by a suitable superset whose attainable sets are
known to be convex, and then determine tight polyhedral
over-approximations, i.e., collections of supporting half-spaces, of
the latter. This yields abstractions characterized in terms of linear
programs. As it turns out, each half-space can be obtained as a
solution of a system of linear equations derived from
\ref{e:TimeDiscreteAutonomousControlSystem:G}
and
of differential equations derived from
\ref{e:TimeContinuousAutonomousControlSystem:F}, respectively.

\subsection{State space quantization}
\label{ss:computation:StateSpaceQuantization}

Quantization of the state space of
\ref{e:TimeDiscreteAutonomousControlSystem:G} is realized by
supplementing \ref{e:TimeDiscreteAutonomousControlSystem:G} with a
quantizer; see Section
\ref{s:intro}.

The system $C$ of quantizer cells is chosen as
follows. We first define a region $K$ of the state space $X$ of
\ref{e:TimeDiscreteAutonomousControlSystem:G} whose local dynamics is
deemed an essential part of the behavior of
\ref{e:TimeDiscreteAutonomousControlSystem:G}, then choose a finite
covering $C'$ of $K$, and finally supplement $C'$ by additional cells
in order to obtain a covering $C$ of $\mathbb{R}^n$.
Intuitively, $K$ is the intended operating range of the quantizer,
whereas cells in $C\setminus C'$ represent overflow symbols.

Of course, what is considered essential local dynamics depends on the
purpose of our
analysis of
\ref{e:TimeDiscreteAutonomousControlSystem:G}, and our choice of $C'$
will also be influenced by other particularities of the problem at
hand; hence a general rule for the choice of the quantizer cannot be
given.
However, the following hypothesis should be fulfilled in order to
ensure the correctness of the algorithm for the computation of
abstractions we are going to present.

\begin{hypothesis}
\label{h:quantizer}
The input alphabet $U$ is finite, and
$C$ is a finite covering of $\mathbb{R}^n$ whose elements are
nonempty convex polyhedra.
For each cell $\Delta \in C' \subseteq C$ there is a superset, denoted
$\widehat \Delta$ in the sequel,
for which $\widehat \Delta \subseteq X$ and attainable sets
$\psi(k,\widehat \Delta,u)$ are convex for
all $k \in \intcc{1;N}$ and all $u \colon \intco{0;k} \to U$,
where $\psi$ denotes the
general solution of \ref{e:TimeDiscreteAutonomousControlSystem:G}, and
$N$, the memory span of the abstraction we seek to obtain.
\end{hypothesis}

In the present section, the above condition plays the role of an
assumption. The non-trivial question of how to verify it
is postponed to Section \ref{s:convex}.

\subsection{Smallest discrete abstractions}
\label{ss:computation:SmallestDiscreteAbstractions}

Let $B \subseteq (U \times C)^{\mathbb{Z}_+}$ denote the behavior of
the quantized system
\ref{e:TimeDiscreteAutonomousControlSystem:G},\ref{e:TimeDiscreteAutonomousControlSystem:Q},
and let $N \in \mathbb{Z}_+$ be given.
The \begriff{$N$-complete hull} $B_N$ of $B$ is the intersection of
all $N$-complete behaviors $B' \subseteq (U \times C)^{\mathbb{Z}_+}$
that contain $B$ as a subset.
Under the name \begriff{strongest $N$-complete approximations},
$N$-complete hulls have been introduced and investigated
by \person{Moor} and his collaborators,
e.g. \cite{Moor99,MoorRaisch99}. It has been shown that
$B \subseteq B_{N+1} \subseteq B_N$ and that $N$-complete hulls are
indeed $N$-complete.%
Thus, the map that assigns
to $B$ its $N$-complete hull $B_N$ is a closure operator
\cite{Schechter97}, and $B_N$ is the smallest discrete abstraction of
memory span $N$ of $B$.
Moreover, $B_N$ admits the following characterization
\cite{Moor99,MoorRaisch99}.

\begin{proposition}
\label{prop:CharNCompleteHull.GQ.pointwise}
Let $C$ be a covering of $\mathbb{R}^n$,
$N \in \mathbb{Z}_+$, $u \colon \intcc{0;N} \to U$,
$\Delta \colon \intcc{0;N} \to C$, $\psi$ the general solution of
\ref{e:TimeDiscreteAutonomousControlSystem:G},
$B$ the behavior of the quantized system
\ref{e:TimeDiscreteAutonomousControlSystem:G},\ref{e:TimeDiscreteAutonomousControlSystem:Q},
and $B_N$ the $N$-complete hull of $B$. Then
\[
B_N
=
\Menge{
w \colon \mathbb{Z}_+ \to W
}{
\forall_{\tau \in \mathbb{Z}_+}\;
(\sigma^{\tau}w)|_{\intcc{0;N}} \in B|_{\intcc{0;N}}
}.
\]
Moreover, the sets $M_0$, \dots, $M_N$ defined by
\begin{equation}
\label{e:prop:CharNCompleteHull.GQ.pointwise}
\negthinspace
M_k
=
\Menge{
\psi(k,x_0,u)
}{
x_0 \in X,
\forall_{\tau \in \intcc{0;k}}\;
\psi(\tau,x_0,u) \in \Delta_{\tau}
}
\end{equation}
satisfy $M_k = \Delta_k \cap G(M_{k-1},u_{k-1})$ for all
$k \in \intcc{1;N}$, and for all $k \in \intcc{0;N}$ we have
$(u,\Delta) \in B|_{\intcc{0;k}}$
iff $M_k \not= \emptyset$.
\end{proposition}

In view of Proposition \ref{prop:CharNCompleteHull.GQ.pointwise},
computing an exact representation of the $k$-complete hull $B_k$ of
the behavior of
\ref{e:TimeDiscreteAutonomousControlSystem:G},\ref{e:TimeDiscreteAutonomousControlSystem:Q}
would require verifying
\begin{equation}
\label{e:Mcondition}
M_k ( u, \Delta ) \not= \emptyset
\end{equation}
for all choices of sequences $u$ and $\Delta$, where
$M_k ( u, \Delta )$ is defined by the right hand side of
\ref{e:prop:CharNCompleteHull.GQ.pointwise}.
To verify
\ref{e:Mcondition}, in turn, one must check whether there is some initial
point $x_0 \in \Delta_0$ such that the trajectory generated by $x_0$
and $u$ visits $\Delta_1, \dots, \Delta_k$ at times $1,\dots,k$;
see also \ref{fig:CharNCompleteHull.GQ.pointwise}(a).
(In fact, $M_k(u,\Delta)$ consists of the values at time $k$ of the
trajectories that satisfy the latter condition.)
To perform that test is, in general, an extremely difficult problem
which may only be exactly solved in rather special situations.
One therefore aims at efficiently computing discrete abstractions that
conservatively approximate the smallest one, $B_k$, and resort to a
test
\begin{equation}
\label{e:MHcondition}
\widehat M_k (u, \Delta) \not= \emptyset
\end{equation}
for some outer approximation
$\widehat M_k (u, \Delta)$ of $M_k (u, \Delta)$,
e.g. \cite{MoorRaisch02}.
On the one hand, the set $\widehat M(u,\Delta)$ should have a simple
structure in order to allow for efficiently testing condition
\ref{e:MHcondition}. On the other hand, that set should approximate
$M(u,\Delta)$  as accurately as possible, since $B_k$ already is an
over-approximation of the actual behavior of the quantized system
\ref{e:TimeDiscreteAutonomousControlSystem:G},\ref{e:TimeDiscreteAutonomousControlSystem:Q}
and the difference $\widehat M(u,\Delta) \setminus M(u,\Delta)$ will
inevitably lead to additional spurious signals.

The following novel characterization of $B_k$, which is not valid in
the more general setting of \cite{Moor99,MoorRaisch99}, will be
crucial in our determination of suitable candidates for
$\widehat M_k (u, \Delta)$.

\begin{proposition}
\label{prop:CharNCompleteHull.GQ.Ginjective}
Let $C$, $N$, $u$, $\Delta$, $\psi$ and $M_k$ be as in Proposition
\ref{prop:CharNCompleteHull.GQ.pointwise} and assume in addition that
$G(\cdot,u)$ is injective for all $u \in U$. Then
\begin{equation}
\label{e:prop:CharNCompleteHull.GQ.Ginjective}
M_k
=
% \Delta_{k}
% \cap
% \bigcap_{\tau=1}^k
% \psi(\tau,\Delta_{k - \tau},u|_{\intco{k - \tau;k}})
\bigcap_{\tau=0}^k
\psi(\tau,X \cap \Delta_{k - \tau},u|_{\intco{k - \tau;k}})
\end{equation}
for all $k \in \intcc{0;N}$.
\end{proposition}

\begin{proof}
\ref{e:prop:CharNCompleteHull.GQ.Ginjective} obviously holds for
$k \in \{0,1\}$, so assume \ref{e:prop:CharNCompleteHull.GQ.Ginjective} holds
for some $k \in \intco{1;N}$. Then, using Proposition
\ref{prop:CharNCompleteHull.GQ.pointwise}, we obtain
\def\myequation{M_{k+1}%
\ifCLASSOPTIONonecolumn\relax\else&\fi%
=\Delta_{k+1} \cap G(M_k,u_k)%
\ifCLASSOPTIONonecolumn\relax\else\\&\fi%
=\Delta_{k+1} \cap%
G\left(%
\bigcap_{\tau=0}^k%
\psi(\tau,X \cap \Delta_{k - \tau},u|_{\intco{k - \tau;k}})%
,u_k\right).%
}
\ifCLASSOPTIONonecolumn%
\[
\myequation
\]
\else
\begin{align*}
\myequation
\end{align*}
\fi
Injectivity of $G(\cdot,u_k)$ implies
$G(A \cap B,u_k) = G(A,u_k) \cap G(B,u_k)$ for any sets $A$ and $B$.
This together with
$G(\psi(\tau,\cdot,u|_{\intco{k - \tau;k}}),u_k) =
\psi(\tau+1,\cdot,u|_{\intcc{k - \tau;k}})$ gives
{\def\minalignsep{0pt}%
\begin{align}
\notag
\negthinspace\negthinspace\negthinspace
M_{k+1}
&=
\Delta_{k+1} \cap
\bigcap_{\tau=0}^k
\psi(\tau+1,X \cap \Delta_{k - \tau},u|_{\intco{k - \tau;k+1}})\\
%\notag
&=
\Delta_{k+1} \cap
\bigcap_{\tau=1}^{k+1}
\psi(\tau,X \cap \Delta_{k + 1 - \tau},u|_{\intco{k + 1 - \tau;k+1}}).
% \\
% &=
% \bigcap_{\tau=0}^{k+1}
% \psi(\tau,X \cap \Delta_{k + 1 - \tau},u|_{\intco{k + 1 - \tau;k+1}}).
\negthinspace\tag*{\qedhere}
\end{align}}
\end{proof}

\subsection{Polyhedral over-approximations of attainable sets}
\label{ss:computation:PolyhedralOverApproximations}

We endow the space $\mathbb{R}^n$ with the standard Euclidean
product $\innerProd{\cdot}{\cdot}$, i.e.,
$
\innerProd{x}{y} = \sum_{i=1}^n x_i y_i
$
for any $x, y \in \mathbb{R}^n$. The derivative and the inverse of a
map $f$ is denoted by $f'$ and $f^{-1}$, respectively, and $f^{\ast}$
is the transpose of $f$ if $f \colon \mathbb{R}^n \to \mathbb{R}^m$ is
linear.

\begin{definition}
\label{def:ComplExt}
For any $C^1$-diffeomorphism $\Phi \colon V \to W$ between open sets
$V,W \subseteq \mathbb{R}^n$, the
\begriff{complementary extension}
$\Phi^{\lozenge} \colon V \times \mathbb{R}^n \to W \times \mathbb{R}^n$
of $\Phi$ is defined
by
\[
\Phi^{\lozenge}(p,v)
=
\left(\Phi(p),\left(\Phi'(p)^{-1}\right)^{\ast}v\right).
\]
\end{definition}
We further define
% $P \colon \mathbb{R}^n \times \mathbb{R}^n \to \mathcal{P}(\mathbb{R}^n)$
% by
\begin{equation}
\label{e:HalfSpace}
P(p,v)
=
\Menge{x \in \mathbb{R}^n}{\innerProd{v}{x-p} \leq 0}
\end{equation}
for all $p, v \in \mathbb{R}^n$ and set
\begin{equation}
\label{e:IntersectionOfHalfSpaces}
P(\Sigma) = \bigcap_{(p,v) \in \Sigma} P(p,v)
\end{equation}
for $\Sigma \subseteq \mathbb{R}^n \times \mathbb{R}^n$.
In words, \ref{e:IntersectionOfHalfSpaces} is the intersection of
the half-spaces \ref{e:HalfSpace} represented by pairs
$(p,v) \in \Sigma$.
\begin{definition}
\label{def:PolyApprox}
A vector $v \in \mathbb{R}^n$ is
\begriff{normal} to $\Omega \subseteq \mathbb{R}^n$ at a boundary
point $p$ of $\Omega$ if $\innerProd{v}{x-p} \leq 0$ for all
$x \in \Omega$.
We call $\Sigma$ an \begriff{outer convex approximation} of
$\Omega$ if $\Omega \subseteq P(\Sigma)$, and
a \begriff{supporting convex approximation} of
$\Omega$, if $p \in \Omega$ and $v$ is normal
to $\Omega$ at $p$, for all $(p,v) \in \Sigma$. A finite outer
(supporting, resp.) convex approximation is \begriff{polyhedral}.
\end{definition}

Let us now return to the problem of suitable candidates
$\widehat M_k (u, \Delta)$ for the test \ref{e:MHcondition}.
If hypothesis \ref{h:quantizer} holds and
$\Delta_0,\dots,\Delta_N \in C'$, we could
define
$\widehat M_k (u, \Delta)$ to be
\begin{equation}
\label{e:MHdef1}
\widehat \Delta_{k}
\cap
\bigcap_{\tau=1}^k
\psi(\tau,\widehat \Delta_{k - \tau},u|_{\intco{k - \tau;k}}),
\end{equation}
and since $\widehat \Delta_{N}$ and all the sets
$\psi(\tau,\widehat \Delta_{N - \tau},u|_{\intco{N - \tau,N}})$ are
convex by \ref{h:quantizer},
the test \ref{e:MHcondition} would be a convex program.
The strategy we actually pursue is to take some
suitable outer polyhedral approximation of \ref{e:MHdef1} for 
$\widehat M_k (u, \Delta)$. Then the convex program
\ref{e:MHcondition} becomes linear, and the sets
$\widehat M_k (u, \Delta)$ enjoy a recursive description.
% , which will
% allow us to efficiently compute and maintain them.

\begin{proposition}
\label{prop:OuterPolyApp}
Assume \ref{h:quantizer} for some $N \in \mathbb{Z}_+$,
as well as \ref{h:GdiffeomorphismC1}, and let
$\Delta \colon \intcc{0;N} \to C'$,
$u \colon \intco{0;N} \to U$ and
$\Sigma, S \colon \intcc{0;N} \to \mathcal{P}(\mathbb{R}^n \times \mathbb{R}^n)$,
where $\mathcal{P}(\cdot)$ denotes the power set.
Assume further that $\Sigma_k$ is a
supporting convex approximation of $\widehat \Delta_k$
for all $k \in \intcc{0;N}$, and
\begin{align}
\label{e:prop:OuterPolyApp:S0}
S_0 &= \Sigma_0,\\
\label{e:prop:OuterPolyApp:Sk}
S_k &= \Sigma_k \cup G(\cdot, u_{k-1})^{\lozenge}(S_{k-1})
\text{\ \ for $k \in \intcc{1;N}$}.
\end{align}
Then, for all $k \in \intcc{1;N}$,
$G(\cdot, u_{k-1})^{\lozenge}(S_{k-1})$ is an outer convex
approximation of
\ifCLASSOPTIONonecolumn\linebreak\fi
$\bigcap_{\tau = 1}^k
\psi(\tau,\widehat \Delta_{k - \tau},u|_{\intco{k - \tau;k}})$, and in
particular, $S_k$ is one of \ref{e:MHdef1}.
\end{proposition}

The above result will enable us to iteratively and efficiently
compute the sets $\widehat M_k (u, \Delta)$ defined earlier. For given
$u$ and $\Delta$, such sets correspond to $P(S_k)$, i.e., to the
intersection of the half-spaces represented by pairs $(p,v) \in S_k$
of points $p$ and normals $v$. In view of this implicit representation
of polyhedra, \ref{e:prop:OuterPolyApp:Sk} says that
$\widehat M_k (u, \Delta)$ is the intersection, and not the union, of
polyhedra $P(\Sigma_k)$ and $P(G(\cdot, u_{k-1})^{\lozenge}(S_{k-1}))$.
Moreover, in contrast to $M_k(u,\Delta)$, the set
$\widehat M_k(u,\Delta)$ is not the intersection of attainable sets of
quantizer cells under the dynamics of
\ref{e:TimeDiscreteAutonomousControlSystem:G}, which is why
Propositions \ref{prop:CharNCompleteHull.GQ.pointwise} and
\ref{prop:CharNCompleteHull.GQ.Ginjective}
cannot be applied to obtain the recursive description in Proposition
\ref{prop:OuterPolyApp}.%

To prove Proposition \ref{prop:OuterPolyApp}
we need the following auxiliary result.

\begin{lemma}
\label{lem:Normal}
Assume $\Phi \colon V \to W$ is a $C^1$-diffeomorphism between open
sets $V,W \subseteq \mathbb{R}^n$ such that both $\Omega \subseteq V$
and $\Phi(\Omega)$ are convex, and let $p \in \Omega$.\\
Then $v \in \mathbb{R}^n$ is normal to $\Omega$ at $p$ iff
$\left(\Phi'(p)^{-1}\right)^{\ast}v$ is normal to $\Phi(\Omega)$ at $\Phi(p)$.
In particular, $\Sigma \subseteq \mathbb{R}^n \times \mathbb{R}^n$ is
a supporting convex approximation of $\Omega$ iff
$\Phi^{\lozenge}(\Sigma)$ is one of $\Phi(\Omega)$.
\end{lemma}

\begin{proof}
Let $v$ be normal to $\Omega$ at $p$ and define
$w = \left(\Phi'(p)^{-1}\right)^{\ast}v$ and
$q = \Phi(p)$
as well as
$\gamma \colon \intcc{0,1} \to \mathbb{R}^n 
\colon t \mapsto \Phi^{-1}(q + t (y-q))$
for some $y \in \Phi(\Omega)$. The map $\gamma$ is
well-defined, differentiable and takes its values in $\Omega$ since
$q,y \in \Phi(\Omega)$ and $\Phi(\Omega)$ is convex. This implies the map
$\alpha$ defined by $\alpha(t) = \innerProd{v}{\gamma(t)-p}$ is
non-positive as $v$ is normal to $\Omega$ at $p$. Furthermore, $\alpha$
is differentiable with $\alpha(0) = 0$, hence
%$0 \geq \alpha'(0) =
%\innerProd{\left(\Phi'(p)^{-1}\right)^{\ast}v}{y-q}$.
$0 \geq \alpha'(0) =
\innerProd{v}{\left( \Phi^{-1} \right)'(q)(y-q)} =
\innerProd{w}{y-q}$.
As $y$ is an arbitrary element of $\Phi(\Omega)$,
$w$ is normal to $\Phi(\Omega)$ at $q$.

For the converse assume $w$ is normal to $\Phi ( \Omega )$ at $q$ and
observe
$\left( ( ( \Phi^{-1} )' ( q ) )^{-1}
\right)^{\ast} w = v$. The first part of this proof applied to
$\Phi^{-1}$ then shows that $v$ is normal to $\Omega$ at $p$.
\end{proof}

\begin{proof}[\IEEEproofname{} of Proposition \ref{prop:OuterPolyApp}]
Let $v \colon \mathbb{Z}_+ \to U$ and observe that $\psi(0,p,v) = p$ and
$\psi(k+1,p,v) = G(\psi(k,p,v),v_k)$ for all $k \in \mathbb{Z}_+$ and
all $p \in X$. Then, by induction, $\psi(k,\cdot,v)$ is a
$C^1$-diffeomorphism between $X$ and an open subset of $X$, which
has two implications. First, by \ref{h:quantizer} and Lemma \ref{lem:Normal},
$\psi(\tau,\cdot,u|_{\intco{k-\tau;k}})^{\lozenge}(\Sigma_{k-\tau})$
is a supporting convex approximation of
$\psi(\tau,\widehat \Delta_{k-\tau},u|_{\intco{k-\tau;k}})$ for all
$k \in \intcc{0;N}$ and all $\tau \in \intcc{0;k}$.
%, thus $\widehat M_k$ contains \ref{e:MHdef1} as a subset.
Second,
% by Lemma \ref{lem:circComplemExt},
we have
\begin{equation}
\label{e:prop:AdjEquDiscrete}
\psi(k+1,\cdot,v)^{\lozenge}
=
G(\cdot,v_k)^{\lozenge}
\circ
\psi(k,\cdot,v)^{\lozenge}
\end{equation}
for all $k \in \mathbb{Z}_+$, where $\circ$ denotes composition
\cite{Schechter97}, since obviously
$(\Phi \circ \Psi)^{\lozenge} = \Phi^{\lozenge} \circ \Psi^{\lozenge}$
for any diffeomorphisms $\Phi$ and $\Psi$ whose
composition $\Phi \circ \Psi$ is well-defined.
We now show
\begin{equation}
\label{e:SkExplizit}
G(\cdot,u_{k-1})^{\lozenge}(S_{k-1})
=
\bigcup_{\tau = 1}^k
\psi(\tau,\cdot,u|_{\intco{k - \tau;k}})^{\lozenge}(\Sigma_{k - \tau})
\end{equation}
for all $k \in \intcc{1;N}$, which proves the proposition. Observe
first that \ref{e:prop:OuterPolyApp:S0} implies
\ref{e:SkExplizit} for $k = 1$, and then assume
\ref{e:SkExplizit} holds for some $k \in \intco{1;N}$. Then
$G(\cdot,u_k)^{\lozenge}(S_k)$ equals
\[
G(\cdot,u_k)^{\lozenge}(\Sigma_k)
\cup
\bigcup_{\tau = 1}^k
G(\cdot,u_k)^{\lozenge}(\psi(\tau,\cdot,u|_{\intco{k - \tau;k}})^{\lozenge}(\Sigma_{k - \tau}))
\]
by \ref{e:prop:OuterPolyApp:Sk}.
The first set in this union is
$\psi(1,\cdot,u|_{\intco{k;k+1}})^{\lozenge}(\Sigma_k)$, and by
\ref{e:prop:AdjEquDiscrete}, the union of the last $k$ sets equals
%\[
$
\bigcup_{\tau = 1}^k
\psi(\tau + 1,\cdot,u|_{\intco{k - \tau;k + 1}})^{\lozenge}(\Sigma_{k - \tau})
$.
%\]
This implies \ref{e:SkExplizit} with $k$ replaced by $k+1$.
\end{proof}

\subsection{Algorithmic solution}
\label{ss:computation:AlgorithmicSolution}

\begin{figure}[t]
\centering
\psfrag{O}[][]{$\Delta$}
\psfrag{imO}[][]{$G(\Delta,u)$}
\psfrag{G}[b][b]{$G(\cdot,u)$}
\psfrag{Go}[b][b]{$G(\cdot,u)^{\lozenge}$}
\psfrag{hO}[bl][bl]{$\widehat\Delta$}
\psfrag{GhO}[br][br]{\footnotesize$G(\widehat\Delta,u)$}
\psfrag{PShO}[b][b]{$P(\Sigma(\widehat\Delta))$}
\psfrag{PGShO}[br][br]{\footnotesize$P(G(\cdot,u)^{\lozenge}(\Sigma(\widehat\Delta)))$}
\includegraphics[width=\ifCLASSOPTIONonecolumn3.5in\else.99\linewidth\fi]{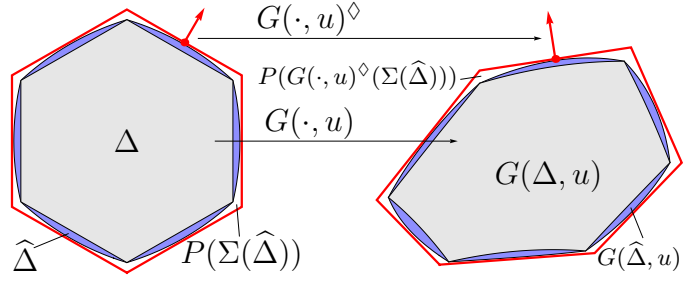}
\caption{\label{fig:HexagonAndStronglyConvexHullAndTheirImage}
Approximation principle underlying the algorithm in
\ref{alg:DiscreteAbstractions}:
Let $u \in U$ and consider a cell
$\Delta \in C'$ whose image
$G(\Delta,u)$ under the nonlinear map $G(\cdot,u)$ is non-convex.
$\Delta$ is conservatively approximated by some set $\widehat \Delta$,
and $\widehat \Delta$ in turn, by a supporting convex polyhedron
$P(\Sigma(\widehat \Delta))$.
$\Sigma(\widehat \Delta)$ is mapped under the
complementary extension $G(\cdot,u)^{\lozenge}$. The result, in
contrast to images under $G(\cdot,u)$, again represents a
convex polyhedron,
$P(G(\cdot,u)^{\lozenge}(\Sigma(\widehat\Delta)))$. By hypothesis
\ref{h:quantizer}, $G(\widehat \Delta,u)$ is convex, which guarantees
$G(\Delta,u) \subseteq G(\widehat \Delta,u) \subseteq
P(G(\cdot,u)^{\lozenge}(\Sigma(\widehat\Delta)))$. The latter
set is the one that is actually computed.
}
\end{figure}

We now present an algorithm for efficiently computing
discrete abstractions for the quantized system
\ref{e:TimeDiscreteAutonomousControlSystem:G},\ref{e:TimeDiscreteAutonomousControlSystem:Q},
which is based on the geometric idea behind Prop.~\ref{prop:OuterPolyApp}.
See also \ref{fig:HexagonAndStronglyConvexHullAndTheirImage}.
Since we are now investigating behaviors, which are sets of signals
$(u,\Delta)$, the analogs of the sets $S_k$ introduced in Proposition
\ref{prop:OuterPolyApp} will depend on $u$ and $\Delta$, which is why
we denote them by
$S(u|_{\intco{\tau;\tau+k}},\Delta|_{\intcc{\tau;\tau+k}})$ in what
follows.

\begin{theorem}
\label{th:AlgorithmDiscreteAbstractions}
Let be $n \in \mathbb{N}$ and $N \in \mathbb{Z}_+$, assume
\ref{h:GdiffeomorphismC1} and \ref{h:quantizer} hold, let
$\Sigma(\widehat \Delta)$ be a supporting polyhedral approximation of
$\widehat \Delta$, for all $\Delta \in C'$, and let the map
$S \colon \bigcup_{k=0}^N
U^{\intco{0;k}} \times C^{\intcc{0;k}}
\to
\mathcal{P}(\mathbb{R}^n \times \mathbb{R}^n)$
be defined by the output of the algorithm in
\ref{alg:DiscreteAbstractions}.
Then, for all $k \in \intcc{0;N}$, the set
\[
\Menge{
(u,\Delta) \in (U \times C)^{\mathbb{Z}_+}
}{
\forall_{\tau \in \mathbb{Z}_+}
P(S(u|_{\intco{\tau;\tau+k}},\Delta|_{\intcc{\tau;\tau+k}}))
\not=\emptyset
},
\]
which is denoted  $B_k$,
is a $k$-complete discrete abstraction of the behavior of
the quantized system
\ref{e:TimeDiscreteAutonomousControlSystem:G},\ref{e:TimeDiscreteAutonomousControlSystem:Q}.
\end{theorem}

\begin{figure}[!t]
{%
\algorithmicindent.6em%
\newcommand{\INPUT}{\item[\textbf{Input:}]}%
\newcommand{\OUTPUT}{\item[\textbf{Output:}]}%
\begin{algorithmic}[1]
\INPUT{$n$, $N$, $U$, $G$, $C$, $C'$;
%$\widehat \Delta$,
$\Sigma(\widehat \Delta)$ for each $\Delta \in C'$.}
%\STMT (* Initialization: *)
\FORALL {$\Delta \in C$}
\label{alg:DiscreteAbstractions:StildeNEW:init1}
 \STATE
 {$S(\Delta) \defas
\begin{cases}
\Sigma(\widehat \Delta), & \text{if $\Delta \in C'$},\\
\emptyset, & \text{otherwise}
\end{cases}
$}
\label{alg:DiscreteAbstractions:StildeNEW:init2}
\ENDFOR
\label{alg:DiscreteAbstractions:StildeNEW:init3}
\FOR{$k=0,\dots,N-1$}
\FORALL {$u \colon \intco{0;k+1} \to U$, $\Delta \colon \intcc{0;k+1} \to C$}
\IF {$S(u|_{\intco{0;k}},\Delta|_{\intcc{0;k}}) = \emptyset$}
\label{alg:DiscreteAbstractions:StildeNEW:ifempty1:if}
\STATE
{$S(u,\Delta)\defas\emptyset$}
\label{alg:DiscreteAbstractions:StildeNEW:ifempty1:then}
\ELSIF
{$\Delta_{k+1} \cap P(G(\cdot,u_{k})^{\lozenge}(
S(u|_{\intco{0;k}},\Delta|_{\intcc{0;k}})
))
 =\emptyset$}
\label{alg:DiscreteAbstractions:TestForEmptyness:if}
\STATE
{$S(u,\Delta)\defas\mathbb{R}^n \times \mathbb{R}^n$}
\label{alg:DiscreteAbstractions:TestForEmptyness:then}
\ELSIF
{$\Delta_{k+1} \notin C'$}
\label{alg:DiscreteAbstractions:TestForOmegaInCprime:if}
\STATE
{$S(u,\Delta)\defas\emptyset$}
\label{alg:DiscreteAbstractions:TestForOmegaInCprime:then}
\ELSE
\STATE
\label{alg:DiscreteAbstractions:StildeNEW}
{$S(u,\Delta)
\defas
\Sigma(\widehat \Delta_{k+1}) \cup
G(\cdot,u_{k})^{\lozenge}(
S(u|_{\intco{0;k}},\Delta|_{\intcc{0;k}})
)$}
\ENDIF
%\STATE
%{$S(\Delta,u)\defas Z$}
\ENDFOR
\ENDFOR
\OUTPUT{$S$}
\end{algorithmic}%
}
\caption{\label{alg:DiscreteAbstractions}Algorithm for the computation
of outer polyhedral approximations of attainable sets that define a discrete abstraction of quantized system
\ref{e:TimeDiscreteAutonomousControlSystem:G},\ref{e:TimeDiscreteAutonomousControlSystem:Q}.}
\end{figure}

\begin{proof}
First observe the operations to be performed by the algorithm in
\ref{alg:DiscreteAbstractions} are well-defined. In particular,
$G(\cdot,u_{k})^{\lozenge}$ on lines
\ref{alg:DiscreteAbstractions:TestForEmptyness:if} and
\ref{alg:DiscreteAbstractions:StildeNEW} of the algorithm shown in
\ref{alg:DiscreteAbstractions} is well-defined by hypotheses
\ref{h:GdiffeomorphismC1} and \ref{h:quantizer}, and
$\Sigma(\widehat \Delta_{k+1})$ on line
\ref{alg:DiscreteAbstractions:StildeNEW} is as well, by the test on
line \ref{alg:DiscreteAbstractions:TestForOmegaInCprime:if}.

Denote the behavior of quantized system
\ref{e:TimeDiscreteAutonomousControlSystem:G},\ref{e:TimeDiscreteAutonomousControlSystem:Q}
by $B$. In
order to show $B_k$ is $k$-complete, let
$u \colon \mathbb{Z}_+ \to U$ and
$\Delta \colon \mathbb{Z}_+ \to C$, and assume that for all
$\tau \in \mathbb{Z}_+$ there is some $(v,\Gamma) \in B_k$ such that
$(\sigma^{\tau} u)|_{\intcc{0;k}} = v|_{\intcc{0;k}}$ and
$(\sigma^{\tau} \Delta)|_{\intcc{0;k}} = \Gamma|_{\intcc{0;k}}$.
This implies
\[
\emptyset
\not=
P(S(v|_{\intco{0;k}},\Gamma|_{\intcc{0;k}}))
=
P(S(u|_{\intco{\tau;\tau+k}},\Delta|_{\intcc{\tau;\tau+k}})),
\]
thus $(u,\Delta) \in B_k$.

Given arbitrary $K \in \intco{0;N}$,
$u \colon \intco{0;K+1} \to U$ and
$\Delta \colon \intcc{0;K+1} \to C$, we now show
\begin{equation}
\label{e:th:AlgorithmDiscreteAbstractions:proof:stern}
% P(S(u|_{\intco{0;k}}, \Delta|_{\intcc{0;k}})) = \emptyset
% \;\;\implies\;\;
% M_{k} = \emptyset
P(S(u|_{\intco{0;k}}, \Delta|_{\intcc{0;k}})) = \emptyset
\;\;\text{implies}\;\;
M_{k} = \emptyset
\end{equation}
for all $k \in \intcc{0;K+1}$, where $M_{k}$ is defined by
\ref{e:prop:CharNCompleteHull.GQ.Ginjective}. According to Propositions
\ref{prop:CharNCompleteHull.GQ.pointwise} and
\ref{prop:CharNCompleteHull.GQ.Ginjective}, this implies
$B \subseteq B_{k}$, and hence, proves the theorem.

From the initialization of $S$ on lines
\ref{alg:DiscreteAbstractions:StildeNEW:init1}-\ref{alg:DiscreteAbstractions:StildeNEW:init3}
of the algorithm and hypothesis \ref{h:quantizer} it follows that
$P(S(\Delta_0)) \not= \emptyset$, thus
\ref{e:th:AlgorithmDiscreteAbstractions:proof:stern} holds for
$k = 0$.
Assume now that \ref{e:th:AlgorithmDiscreteAbstractions:proof:stern} holds
for all $k \in \intcc{0;K}$, for some $K \in \intco{0;N}$, as
well as $P(S(u,\Delta)) = \emptyset$.
In view of lines \ref{alg:DiscreteAbstractions:StildeNEW:ifempty1:if}
and \ref{alg:DiscreteAbstractions:StildeNEW:ifempty1:then}, this
implies
$S(u|_{\intco{0;k}},\Delta|_{\intcc{0;k}}) \not= \emptyset$
for all $k \in \intcc{0;K+1}$. We further obtain
\begin{equation}
\label{e:th:AlgorithmDiscreteAbstractions:proof:TestForEmptyness:if}
\Delta_{k + 1}
\cap
P(G(\cdot,u_k)^{\lozenge}
(S(u|_{\intco{0;k}},\Delta|_{\intcc{0;k}})))
=
\emptyset
\end{equation}
for $k = K$. Otherwise, $S(u,\Delta)$ would have been assigned
its value on line \ref{alg:DiscreteAbstractions:StildeNEW}, i.e.,
{%
\def\myequation{\label{e:th:AlgorithmDiscreteAbstractions:proof:StildeNEW}
S(u|_{\intco{0;k+1}},\Delta|_{\intcc{0;k+1}})
=
\ifCLASSOPTIONonecolumn\relax\else\\\fi
\Sigma(\widehat \Delta_{k+1})
\cup
G(\cdot,u_k)^{\lozenge}
(S(u|_{\intco{0;k}},\Delta|_{\intcc{0;k}}))
}
\ifCLASSOPTIONonecolumn%
\begin{equation}
\myequation
\end{equation}
\else
\begin{multline}
\myequation
\end{multline}%
\fi%
}%
for $k = K$; thus the left hand side of
\ref{e:th:AlgorithmDiscreteAbstractions:proof:TestForEmptyness:if} would
be a subset of $P(S(u,\Delta)) = \emptyset$, which is a contradiction.

Now assume
\ref{e:th:AlgorithmDiscreteAbstractions:proof:TestForEmptyness:if} also
holds for some $k \in \intco{0;K}$.
Then
$S(u|_{\intco{0;k+1}},\Delta|_{\intcc{0;k+1}})$ is assigned
its value on line
\ref{alg:DiscreteAbstractions:TestForEmptyness:then}, hence
$P(S(u|_{\intco{0;k+1}},\Delta|_{\intcc{0;k+1}})) = \emptyset$.
From this and \ref{e:th:AlgorithmDiscreteAbstractions:proof:stern} we
obtain $M_{k+1} = \emptyset$, and Proposition
\ref{prop:CharNCompleteHull.GQ.pointwise} yields
$M_{K+1} = \emptyset$. Thus,
\ref{e:th:AlgorithmDiscreteAbstractions:proof:stern} holds for
$k = K+1$.

If, on the other hand,
\ref{e:th:AlgorithmDiscreteAbstractions:proof:TestForEmptyness:if}
does not hold for any $k \in \intco{0;K}$, then, in view of lines
\ref{alg:DiscreteAbstractions:StildeNEW:init2},
\ref{alg:DiscreteAbstractions:StildeNEW:ifempty1:if},
\ref{alg:DiscreteAbstractions:StildeNEW:ifempty1:then},
\ref{alg:DiscreteAbstractions:TestForOmegaInCprime:if} and
\ref{alg:DiscreteAbstractions:TestForOmegaInCprime:then}, we obtain
$\Delta_{\tau} \in C'$ for all $\tau \in \intcc{0;K}$,
$S(\Delta_0) = \Sigma(\widehat \Delta_0)$, as well as
\ref{e:th:AlgorithmDiscreteAbstractions:proof:StildeNEW} for all
$k \in \intco{0;K}$. Proposition \ref{prop:OuterPolyApp} then
shows that the left hand side of
\ref{e:th:AlgorithmDiscreteAbstractions:proof:TestForEmptyness:if} for
$k = K$ contains $M_{K+1}$ as a subset, so
\ref{e:th:AlgorithmDiscreteAbstractions:proof:stern} holds for
$k = K+1$, and we are done.
\end{proof}

\begin{corollary}
\label{cor:th:AlgorithmDiscreteAbstractions}
Under the hypotheses of Theorem \ref{th:AlgorithmDiscreteAbstractions}
and for all $k \in \intcc{0;N}$, the requirement
\[
\exists_{\Delta \colon \mathbb{Z}_+ \to C}
\left(
(u,\Delta)\in B_k
\text{\ and\ }
\forall_{k \in \mathbb{Z}_+}
x_k \in \Delta_k
\right)
\]
for sequences $(u,x) \colon \mathbb{Z}_+ \to U \times \mathbb{R}^n$
defines an abstraction of the behavior of the discrete-time system
\ref{e:TimeDiscreteAutonomousControlSystem:G}.
\end{corollary}

We remark that the algorithm in \ref{alg:DiscreteAbstractions} contains
just two nontrivial operations which
need to be performed repeatedly, namely, the determination of the set
\begin{equation}
\label{e:Stilde}
G(\cdot,u_{k})^{\lozenge}(S(u|_{\intco{0;k}},\Delta|_{\intcc{0;k}})),
\end{equation}
which appears on lines \ref{alg:DiscreteAbstractions:TestForEmptyness:if} and
\ref{alg:DiscreteAbstractions:StildeNEW}, and the test for emptiness on line
\ref{alg:DiscreteAbstractions:TestForEmptyness:if}.
The latter can be efficiently performed using linear programming
techniques, since both $\Delta_{k+1}$ and
$P(G(\cdot,u_{k})^{\lozenge}(
S(u|_{\intco{0;k}},\Delta|_{\intcc{0;k}})
))$ are convex polyhedra.
According to Definition \ref{def:ComplExt}, the former operation
requires an evaluation of the function $G(\cdot,u_k)$ and the solution
of a linear system of equations for each element
$(p,v) \in S(u|_{\intco{0;k}},\Delta|_{\intcc{0;k}})$,
\begin{subequations}
\label{e:TimeDiscreteAutonomousControlSystem:G:adjoint}
\begin{align}
\label{e:TimeDiscreteAutonomousControlSystem:G:adjoint:a}
\tilde p &= G(p,u_k),\\
\label{e:TimeDiscreteAutonomousControlSystem:G:adjoint:b}
v &= D_1 G(p,u_k)^{\ast} \tilde v,
\end{align}
\end{subequations}
in order to obtain an element $(\tilde p,\tilde v)$ of \ref{e:Stilde},
which represents one half-space in the outer convex approximation
\ref{e:Stilde}. Here, $D_i f$ denotes the partial derivative of $f$
with respect to the $i$th argument.

In order to estimate the computational complexity of the algorithm in
\ref{alg:DiscreteAbstractions}, we assume for simplicity that each of
the sets $\widehat \Delta$ is approximated by $m$ supporting
half-spaces, i.e., $m = | \Sigma ( \widehat \Delta ) |$ for all
$\Delta \in C'$, where $|\cdot|$ denotes cardinality. It is then easy
to see that for any given $k \in \intco{0;N}$, the number of
half-spaces needed to define all the sets \ref{e:Stilde} is at most
$m |C| |U|^{k+1}$, and these sets are represented by at most $(k+1)m$
half-spaces each.
To estimate the number of tests for emptiness,
% on line \ref{alg:DiscreteAbstractions:TestForEmptyness:if}
we additionally assume that there is a constant $\lambda > 0$, independent
of $|C|$, such that the following holds.
For any given set of the form \ref{e:Stilde}, the test on line
\ref{alg:DiscreteAbstractions:TestForEmptyness:if} has to be
performed for at most $\lambda$ cells $\Delta_{k+1} \in C$, and the
set of these candidate cells can be provided in constant time.
This holds, in particular, if the cells in $C'$ are congruent compacta
arranged in a regular grid. Then, for any given $k$, the number of
tests on line \ref{alg:DiscreteAbstractions:TestForEmptyness:if} is
bounded by $\lambda^{k+1} |C| |U|^{k+1}$. Note here that the values
$\emptyset$ and $\mathbb{R}^n \times \mathbb{R}^n$, which may be
assigned on lines \ref{alg:DiscreteAbstractions:TestForEmptyness:then}
and \ref{alg:DiscreteAbstractions:TestForOmegaInCprime:then}, play a
role similar to zeros in sparse matrices, and thus,
these values do not need to be stored and computations on them do not
need to actually be performed \cite{i06Fill}.

In summary, the algorithm in \ref{alg:DiscreteAbstractions} requires
the solution of $O( m | C | | U |^N )$ instances of
\ref{e:TimeDiscreteAutonomousControlSystem:G:adjoint} and of
$O( \lambda^N | C | | U |^N )$ linear feasibility problems in $n$ variables with
at most $(N+1)m$ inequalities each, where $O(\cdot)$ is the usual
asymptotic notation \cite{Gruber07}.
The parameters $m$ and $|C|$ depend on the dimension $n$ of
the state space of \ref{e:TimeDiscreteAutonomousControlSystem:G}, with
$|C|$ typically
growing exponentially. Therefore, the computational effort has to be
expected to grow rapidly with $n$, a problem that is
common to all grid based methods for the computation of abstractions,
e.g.
\cite{Hsu87,Osipenko07,Tabuada08,PolaTabuada09,JaulinWalter97}.

We remark that apart from an increasing
computational effort, application of the algorithm shown in
\ref{alg:DiscreteAbstractions} in dimensions exceeding $2$ does not
pose any particular difficulties. This is obvious for the operations
on lines \ref{alg:DiscreteAbstractions:TestForEmptyness:if} and
\ref{alg:DiscreteAbstractions:StildeNEW}, which we have already
discussed, and also holds for the remaining operations.
Specifically, for the computation of the supporting polyhedral
approximations $\Sigma( \widehat \Delta )$ on line
\ref{alg:DiscreteAbstractions:StildeNEW:init2}, several methods are
available for a large class of sets $\widehat \Delta$ \cite{Gruber07}.

Finally, we would like to emphasize again that convexity of certain
attainable sets is an important requirement for the correctness of the
algorithm in \ref{alg:DiscreteAbstractions}, see hypothesis
\ref{h:quantizer}.
While for linear systems that requirement is always met by the choice
$\widehat \Delta = \Delta$, the results of section \ref{s:convex} will
show how to meet it in the presence of
non-linearities.

\subsection{Sampled systems}
\label{ss:SampledSystems}

Here we consider the case that
\ref{e:TimeDiscreteAutonomousControlSystem:G} arises from a
continuous-time system \ref{e:TimeContinuousAutonomousControlSystem:F}
under sampling.
More formally, let a continuous-time control system
\ref{e:TimeContinuousAutonomousControlSystem:F}
with $F \colon X \times V \to \mathbb{R}^n$ and
a set $U$ of input signals be given, where
$X \subseteq \mathbb{R}^n$,
$V \subseteq \mathbb{R}^m$,
each $u \in U$ is a piecewise continuous map
$u \colon \intcc{0,T} \to V$, and
$T > 0$ is the \begriff{sampling period}.
A map $v \colon \mathbb{R}_+ \to V$
is an \begriff{admissible input signal} for
\ref{e:TimeContinuousAutonomousControlSystem:F}, generated by
$u \colon \mathbb{Z}_+ \to U$, if
$
v(t) = u_k(t - k T)
$
for all $k \in \mathbb{Z}_+$ and all
$t \in \intco{k T,(k + 1) T}$.
The set of admissible input signals for
\ref{e:TimeContinuousAutonomousControlSystem:F} is denoted
$\mathcal{V}$ in the sequel.
We assume the following.
\begin{hypothesis}
\label{h:FstetigC1LoesungenAufRplus}
$X \subseteq \mathbb{R}^n$ is open,
the right hand side $F$ of
\ref{e:TimeContinuousAutonomousControlSystem:F} is continuously
differentiable with respect to its first argument and continuous.
Furthermore,
for any $x_0 \in X$ and any admissible input signal
$v \in \mathcal{V}$, the solution of the initial value problem
composed of \ref{e:TimeContinuousAutonomousControlSystem:F} and the
initial condition
%\begin{equation}
%\label{e:TimeContinuousAutonomousControlSystem:IC}
$x(0)=x_0$
%\end{equation}
is extendable to the entire time axis $\mathbb{R}_{+}$.
\end{hypothesis}

Discrete-time system \ref{e:TimeDiscreteAutonomousControlSystem:G} is called
the \begriff{sampled system} associated with
\ref{e:TimeContinuousAutonomousControlSystem:F} if its right hand side
is given by
\[
G(x,u) = \varphi(T,x,u)
\]
for all $x \in X$ and all $u \in U$, where $\varphi$ is the
\begriff{general solution} of
\ref{e:TimeContinuousAutonomousControlSystem:F}, i.e.,
$\varphi(t,x_0,v)$ is the value at time $t$ of the solution of the
initial value problem composed of
\ref{e:TimeContinuousAutonomousControlSystem:F} and the initial
condition
%,\ref{e:TimeContinuousAutonomousControlSystem:IC}.
$x(0)=x_0$.
% In complete analogy to the discrete-time case, it is not necessary to
% specify $v$ on the whole time axis.

Obviously, the sampled system
\ref{e:TimeDiscreteAutonomousControlSystem:G} associated with
\ref{e:TimeContinuousAutonomousControlSystem:F} fulfills
\ref{h:GdiffeomorphismC1} if
\ref{e:TimeContinuousAutonomousControlSystem:F} fulfills
\ref{h:FstetigC1LoesungenAufRplus}, and
$
\psi(k,x,u) = \varphi(k T,x,v)
$
for all  $x \in X$ and all $k \in \mathbb{Z}_+$, whenever $v$ is an
admissible input signal for
\ref{e:TimeContinuousAutonomousControlSystem:F} generated by the
sequence $u \colon \mathbb{Z}_+ \to U$ and $\psi$ is the general
solution of \ref{e:TimeDiscreteAutonomousControlSystem:G}.
Hence, our results for \ref{e:TimeDiscreteAutonomousControlSystem:G},
including the algorithm in \ref{alg:DiscreteAbstractions}, can be
directly applied to the sampled system
\ref{e:TimeDiscreteAutonomousControlSystem:G} associated with
\ref{e:TimeContinuousAutonomousControlSystem:F} if the latter
satisfies hypothesis \ref{h:FstetigC1LoesungenAufRplus}.
In particular, \ref{e:TimeDiscreteAutonomousControlSystem:G:adjoint}
can be efficiently solved even though the right hand side $G$ of the
sampled system \ref{e:TimeDiscreteAutonomousControlSystem:G}
is not explicitly given. The solution is obtained through solving an
initial value problem in a $2n$-dimensional ordinary differential
equation (ODE) over a single sampling interval:

\begin{proposition}
\label{prop:TimeContinuousAutonomousControlSystem:F}
Let \ref{e:TimeDiscreteAutonomousControlSystem:G} be the sampled
system associated with \ref{e:TimeContinuousAutonomousControlSystem:F}
for sampling period $T > 0$, and
assume \ref{h:FstetigC1LoesungenAufRplus}. Then
\[
G(\cdot,u)^{\lozenge}(p,v)
=
\left(
x(T),
y(T)
\right)
\]
for all $u \in U$, $p \in X $ and $v \in \mathbb{R}^n$, where $(x,y)$
denotes the solution of the initial value problem
\begin{subequations}
\label{e:TimeContinuousAutonomousControlSystem:F:adjoint}
\begin{align}
\label{e:TimeContinuousAutonomousControlSystem:F:adjoint:a}
\dot x(t)
&=
F(x(t),u(t)),\\
\label{e:TimeContinuousAutonomousControlSystem:F:adjoint:b}
\dot y(t)
&=
-D_1F(x(t),u(t))^{\ast} y(t),\\
%\label{e:TimeContinuousAutonomousControlSystem:F:adjoint:c}
x(0) &= p,\;\;
\label{e:TimeContinuousAutonomousControlSystem:F:adjoint:d}
y(0) = v.
\end{align}
\end{subequations}
\end{proposition}

\begin{proof}
Assume $u$ is continuous and let $\varphi$ denote the general solution
of \ref{e:TimeContinuousAutonomousControlSystem:F}.
\ref{h:FstetigC1LoesungenAufRplus} implies $\varphi$ is
continuously differentiable, $X \defas D_2 \varphi(\cdot,p,u)$
fulfills the variational equation
$\dot X(t) = D_1F(x(t),u(t)) X(t)$,
and $X(0) = \id$,
e.g. \cite{Hartman02}. In particular, $G(\cdot,u)$
is a $C^1$-diffeomorphism, so its complementary extension
$G(\cdot,u)^{\lozenge}$ is well-defined.
Furthermore, $\left(X(\cdot)^{-1}\right)^{\ast}$ fulfills the adjoint
equation \ref{e:TimeContinuousAutonomousControlSystem:F:adjoint:b},
e.g. \cite{Hartman02}. Thus,
$\left(D_1 G(p,u)^{-1}\right)^{\ast}v =
\left(X(T)^{-1}\right)^{\ast}v =
y(T)$.
The extension to piecewise continuous $u$ is obvious.
\end{proof}

\section{Convexity of attainable sets}
\label{s:convex}
\begin{figure*}[!t]
\hspace*{\fill}
\begin{minipage}{.61\linewidth}
\noindent
\centering
\psfrag{O}[b][b]{$\Omega$}
\psfrag{x}[t][t]{$x$}
\psfrag{y}[bl][bl]{$x + h + v \cdot \eta_x(h)$}
\psfrag{hv}[bl][bl]{}
\psfrag{h}[br][br]{$h$}
\psfrag{v}[bl][bl]{$v$}
\psfrag{OO}[r][r]{$\varphi(T,\Omega,u)$}
\includegraphics[width=\linewidth]{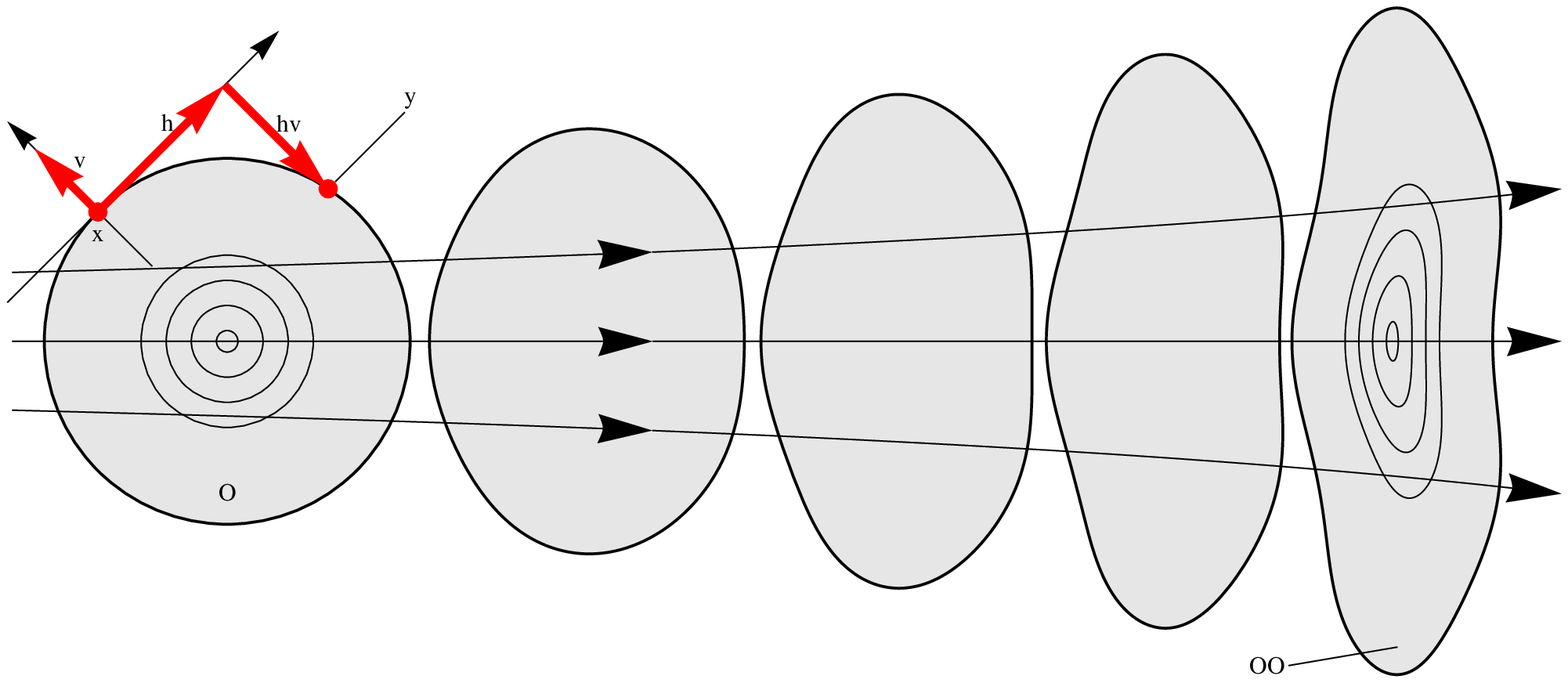}\\
(a)
\end{minipage}
\hspace*{\fill}
\begin{minipage}{.27\linewidth}
\centering
\psfrag{O}[b][b]{$\Delta$}
\psfrag{hO}[][]{\ifCLASSOPTIONonecolumn\small\fi$\widehat\Delta$}
\psfrag{PShO}[br][br]{\ifCLASSOPTIONonecolumn\footnotesize\else\small\fi$P(\Sigma(\widehat\Delta))$}
\psfrag{e}[l][l]{$\varepsilon$}
\psfrag{s}[b][b]{$s$}
\includegraphics[width=\linewidth]{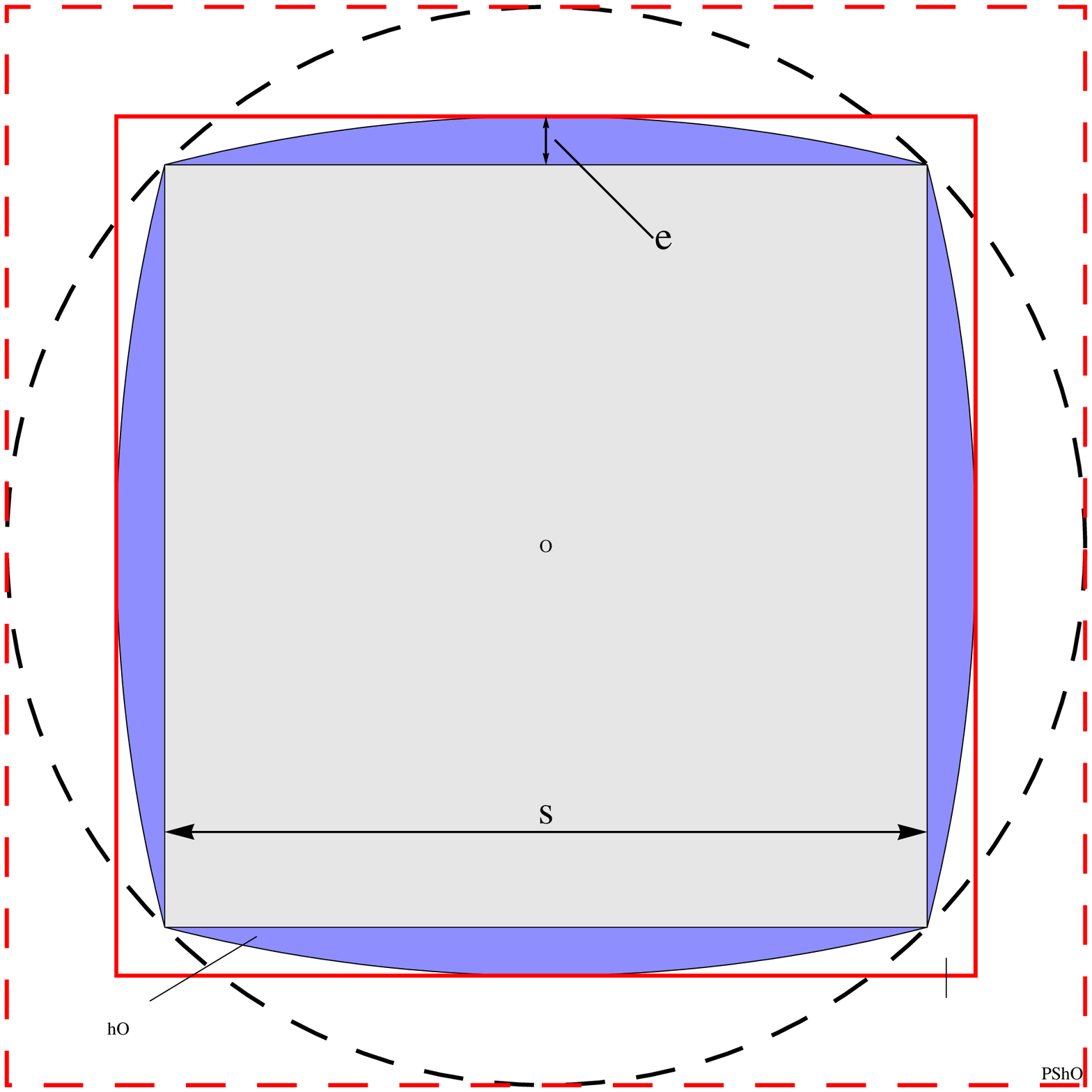}\\
(b)
\end{minipage}
\hspace*{\fill}
\caption{\label{fig:GeometricIntuitionConvexityForSmallTimes}
(a)
Geometric idea behind our derivations in Section \ref{s:convex} for
the continuous-time system
\ref{e:TimeContinuousAutonomousControlSystem:F}:
If $\Omega$ is an Euclidean ball of radius $r$ and the right hand side
$F$ of \ref{e:TimeContinuousAutonomousControlSystem:F} is of class
$C^{1,1}$, then the attainable set $\varphi(T,\Omega,u)$ is convex
whenever $T$ or $r$ is small enough. This idea does not extend to
smooth, strictly convex sets $\Omega$, nor to right hand sides $F$ of
class $C^1$.
(b)
The smallest intersection $\widehat \Delta$ of a finite number of
closed balls of radius $r$ containing a polyhedron $\Delta$ is a
subset of any $r$-convex ellipsoid (dashed) containing $\Delta$.
}
\end{figure*}
In this section we will investigate the convexity of attainable
sets of control systems
\ref{e:TimeDiscreteAutonomousControlSystem:G} and
\ref{e:TimeContinuousAutonomousControlSystem:F}.
To begin with, we briefly explain the geometric idea behind our
derivations using the example of the continuous-time system
\ref{e:TimeContinuousAutonomousControlSystem:F}. Consider a set
$\Omega \subseteq \mathbb{R}^n$ of class $C^2$, denote its boundary by
$\boundary \Omega$, let $x \in \boundary \Omega$ be an arbitrary
boundary point,
and let $v$ be a unit-length normal to $\Omega$ at $x$.
Let $\eta_x$ be a $C^2$-map defined on the tangent space of
$\boundary \Omega$ at $x$ that represents the boundary
$\boundary \Omega$ locally about the point $x$ in local coordinates,
with the origin located at $x$. That is, $\eta_x$ is such that in a
neighborhood of $x$ we have $y \in \boundary \Omega$ iff there is a
tangent vector $h$ to $\boundary \Omega$ at $x$ such that
$y = x + h + v \cdot \eta_x(h)$.
See \ref{fig:GeometricIntuitionConvexityForSmallTimes}(a).
The map $\eta_x$ is locally uniquely determined by
$\Omega$ and satisfies $\eta_x(0) = 0$ and
$\eta_x'(0) = 0$. Moreover, $\Omega$ is convex iff $\eta_x''(0)$ is
negative semi-definite for every $x \in \boundary \Omega$
\cite{i07Convex}, where $\eta_x''$ denotes the second order derivative
of $\eta_x$. Consequently, if the right hand side $F$ of
\ref{e:TimeContinuousAutonomousControlSystem:F} is of class $C^2$, the
issue of convexity of attainable sets can be decided
by studying the evolution of $\eta_x''(0)$ under the
dynamics of \ref{e:TimeContinuousAutonomousControlSystem:F}. The
realization of that idea in the case that $\Omega$ is a Euclidean ball
centered at some point $x_0$
reveals that the
attainable set $\varphi(T,\Omega,u)$ is convex whenever the time
$T$ or the radius $r$ of $\Omega$ is small enough \cite{i07Convex}.
See \ref{fig:GeometricIntuitionConvexityForSmallTimes}(a).
Moreover, bounds on $T$ and $r$ can be derived from properties of the
right hand side $F$ of
\ref{e:TimeContinuousAutonomousControlSystem:F}. These results
generalize to ellipsoids at the place of balls and
can also be generalized to the case of $C^{1,1}$-smoothness by
using suitably generalized second order derivatives \cite{i07Convex}.
They do not, however, extend to smooth, strictly convex sets, nor
to right hand sides of class $C^1$. In particular, attainable sets
$\varphi(T,\Omega,u)$ can then be non-convex for any time $T > 0$.

In the present section, we will replace ellipsoids by strongly convex
sets. Here, the
set $\Omega \subseteq \mathbb{R}^n$ is \begriff{strongly convex of
  radius $r$}, or \begriff{$r$-convex} for short, if it is an
intersection of a family of closed balls of radius $r > 0$, i.e., if
\begin{equation}
\label{e:def:StrongConvexity}
\Omega
=
\bigcap_{x \in M}
\cBall(x,r)
\end{equation}
for some $M \subseteq \mathbb{R}^n$, where $\cBall(x,r)$ denotes the
closed Euclidean ball of radius $r$ centered at $x$.
$\Omega$ is \begriff{strongly convex} if it is $r$-convex for some
$r > 0$ \cite{FrankowskaOlech81,Polovinkin96}.

In view of the method proposed in Section \ref{s:computation}, the use
of strongly convex rather than ellipsoidal supersets of quantizer
cells allows for more precise approximations of attainable sets, and
thus, for more accurate abstractions. Indeed, in
contrast to the case of ellipsoidal supersets, the error $\varepsilon$
by which $P(\Sigma(\widehat \Delta))$ approximates the cell
$\Delta$ in \ref{fig:GeometricIntuitionConvexityForSmallTimes}(b) is
quadratic in the edge length $s$ of $\Delta$, and it can
be shown that this property carries over to approximations of
attainable sets of $\Delta$. See also
\ref{fig:HexagonAndStronglyConvexHullAndTheirImage}.%

We will also extend previous results to the discrete-time
case \ref{e:TimeDiscreteAutonomousControlSystem:G}.
Our results
not only allow verification of the requirements for the correctness of
the algorithm of section \ref{s:computation} when a particular
quantizer is given, but they also help in constructing admissible
quantizers. In particular, we show that if hypothesis
\ref{h:GdiffeomorphismC1} is fulfilled with sufficient smoothness
and $Y \subseteq \mathbb{R}^n$ is a compact, full-dimensional,
convex polyhedron, then choosing
sufficiently small scaled and translated copies of $Y$ as
operating range quantizer cells will guarantee that the stated
requirements of section \ref{s:computation} are met.

As the problems investigated in this section become trivial in
dimension $1$, we assume a multidimensional setting, i.e., $n \geq 2$
throughout this section.

\subsection{Convexity of diffeomorphic images of strongly convex sets}
\label{ss:convex:ConvexityOfDiffeomorphicImages}

We first present a sufficient condition for the diffeomorphic image of
a strongly convex set to be itself convex, which will be the basis of
all subsequent results. In what follows, we write $x \perp y$
if $\innerProd{x}{y} = 0$ and denote by $\| \cdot \|$ the
Euclidean norm of both vectors and linear maps. The
closure, interior and boundary of a set $M \subseteq \mathbb{R}^n$ is
denoted $\closure M$, $\interior M$ and $\boundary M$, respectively.
We set $f h^k \defas f(h,\dots,h)$ if $f$ is $k$-linear.

\begin{proposition}
\label{prop:StrongConvexityDiffeomorph}
Let $\Phi \colon U \to V$ be a $C^{1,1}$-diffeomorphism between open sets
$U, V \subseteq \mathbb{R}^n$
and $\Omega \subseteq U$ be $r$-convex, $\Omega \not= \mathbb{R}^n$.
Assume that for each
$x \in \boundary \Omega$ there is a unit length normal $v$ to $\Omega$
at $x$ such that
\begin{equation}
\label{e:prop:StrongConvexityDiffeomorph:L1}
\limsup_{t \to 0, t > 0}
\frac{\innerProd{v}{\Phi'(x)^{-1}(\Phi'(x+t \xi) - \Phi'(x) ) \xi}}{t}
<
\frac{1}{r} \| \xi \|^2
\end{equation}
holds for all $\xi \perp v$, $\xi \not= 0$.
Then $\Phi(\Omega)$ is convex.
\end{proposition}

It should be noted that the left hand side of
\ref{e:prop:StrongConvexityDiffeomorph:L1} equals
$\innerProd{v}{\Phi'(x)^{-1} \Phi''(x) \xi^2}$ if $\Phi$ is of class $C^2$.
In addition, if \ref{e:prop:StrongConvexityDiffeomorph:L1} were not a
strict inequality and if $\Omega$ is a closed ball of radius $r$, the
condition is known to be necessary and sufficient for the convexity of
$\Phi(\Omega)$ \cite{i07Convex}.
However, this does not prove the proposition,
despite representation \ref{e:def:StrongConvexity}.
Indeed, \ref{e:prop:StrongConvexityDiffeomorph:L1} is only
required to hold at boundary points $x$ of $\Omega$. Hence, even if
all balls that appear in \ref{e:def:StrongConvexity} are in the domain
of definition $U$ of $\Phi$,
\ref{e:prop:StrongConvexityDiffeomorph:L1} may very well be violated
at boundary points of these balls.

In the course of proving Proposition
\ref{prop:StrongConvexityDiffeomorph} we will say that
$\Omega$ is \begriff{weakly supported at
$p \in \boundary \Omega$ locally} whenever there is some neighborhood
$U$ of $p$ and a non-zero normal to $\Omega \cap U$ at $p$
\cite{Valentine64}. Further, we will call $f \colon U \to \mathbb{R}$
a \begriff{$C^{1,1}$-submersion on its zero set} if the following
holds:
$f$ is continuous on the open set $ U \subseteq \mathbb{R}^n$, and for
every zero $x$ of $f$, $f$ is of class $C^{1,1}$ on a neighborhood of
$x$ and $f'(x)$ is surjective.

\begin{lemma}
\label{lem:C11QuadSupport:b}
Let $g \colon U \subseteq \mathbb{R}^n \to \mathbb{R}$ be a
$C^{1,1}$-submersion on its zero set,
$\Omega = \Menge{x \in U}{g(x) \leq 0}$, and assume $g(0) = 0$ and
\begin{equation}
\label{e:lem:C11QuadSupport:b}
\liminf_{t \to 0, t > 0}
 \frac{g'(t h)h}{t}
> 0
\end{equation}
for all $h \in \ker g'(0) \setminus \{ 0\}$. Then $\Omega$ is weakly
supported at $0$ locally.
\end{lemma}

\begin{proof}
Set $v = g'(0)^{\ast} / \| g'(0) \|$.
An application of the implicit function theorem to the equation
$g(h + \lambda v) = 0$ for
$h \in \ker g'(0)$ and $\lambda \in \mathbb{R}$
shows that $\Omega$ can be represented locally about $0$ by a map
$\eta \colon W \subseteq \ker g'(0) \to \mathbb{R}$ of class $C^{1,1}$,
$W$ an open neighborhood of the origin.
That is, for $h$ and $\lambda$ small enough we have
$h + \lambda v \in \Omega$ iff $\lambda \leq \eta(h)$.
Combine this with the identity
$\innerProd{v}{h + \lambda v} = \lambda$ and the definition of weak
local support to see that it suffices to show $\eta(h) \leq 0$ for all
sufficiently small $h$.
In order to prove the latter, first observe that
% $\eta$ is of class $C^{1,1}$ as $g$ is. Moreover,
$\eta(0) = 0$ and $\eta'(0) = 0$. Then differentiate the
identity $g(h + \eta(h) v) = 0$ with respect to $h$ and use the
Lipschitz continuity of $g'$ to obtain
\[
\limsup_{t \to 0, t > 0}
%\frac{\eta'(t h)h}{t}
\eta'(t h)h/t
=
-
%\frac{1}{\| g'(0) \|}
\| g'(0) \|^{-1}
\liminf_{t \to 0, t > 0}
%\frac{g'(t h)h}{t}
g'(t h)h/t
\]
for all $h \in \ker g'(0)$.
If $g$ is of class $C^2$, then so is $\eta$, the left hand side of the
latter equation equals $\eta''(0)h^2$, and the claim follows from
\ref{e:lem:C11QuadSupport:b}. If $g$ is
merely $C^{1,1}$, use again \ref{e:lem:C11QuadSupport:b} and apply
\cite[Theorem 3.2]{BednarikPastor04}.
\end{proof}

\begin{proof}[Proof of Proposition
  \ref{prop:StrongConvexityDiffeomorph}]
The claim is trivial for $\Omega = \emptyset$ and $\Omega$ a
singleton, so we assume $\Omega$ contains at least two points.
Then $\Omega$ has nonempty interior by Lemma
\ref{lem:sConvexUnitBall} in the Appendix. In addition, $\Omega$ is
compact and convex. Hence $\Omega = \closure ( \interior ( \Omega ) )$
and $\interior ( \Omega )$ is connected, and these properties are
preserved under diffeomorphisms.
Moreover,
$\interior ( \Phi ( \Omega ) ) = \Phi ( \interior ( \Omega ) )$ since
$\Phi$ is a diffeomorphism.

We will show below that $\interior ( \Phi ( \Omega ) )$ is weakly
supported at each of its boundary points locally. Then, since that set
is also open and connected, it is convex
% by a  well-known Theorem of \person{Tietze}
\cite[Theorem 4.10]{Valentine64}, which implies its
closure $\Phi ( \Omega )$ is also convex.

Let $x \in \boundary \Omega$ be arbitrary and $v$ be as
in the statement of the
\label{13737A:11:iv}%
proposition,
and assume $x = \Phi(x) = 0$ without loss of generality.
Then $\Omega \subseteq \cBall(-r v, r)$ since $\Omega$ is $r$-convex
and compact \cite[Proposition 3.1]{FrankowskaOlech81}. Now
% let $\circ$ denote the composition of maps \cite{Schechter97} and
define
$f(z) = \| z + r v \|^2 - r^2$ and $g = f \circ \Phi^{-1}$ to observe
that $g(y) \leq 0$ is equivalent to $y \in \Phi ( \cBall(-r v, r) )$,
hence
\begin{equation}
\label{e:prop:StrongConvexityDiffeomorph:proof:1}
\Phi( \Omega )
\subseteq
\Menge{y \in V}{g(y) \leq 0}.
\end{equation}
We claim that the set on the right hand side of
\ref{e:prop:StrongConvexityDiffeomorph:proof:1}
is weakly supported at the origin locally. To prove this, first
observe that $g$ is a $C^{1,1}$-submersion on its zero set since $f$
is one and $\Phi$ is a $C^{1,1}$-diffeomorphism, and that $g(0) = 0$.
Then differentiate the identity $f = g \circ \Phi$, observe
$f'(t \xi) \xi / t = 2 \| \xi \|^2$, and use the Lipschitz continuity
of $g'$ and the continuity of $\Phi'$ to see that $2 \| \xi \|^2$
equals
\[
\liminf_{t \to 0, t > 0}
\left(
%\frac{g'(th)h}{t}
g'(th)h/t
+
2 r
%\frac{\innerProd{v}{\Phi'(0)^{-1}(\Phi'(t \xi) - \Phi'(0)) \xi}}{t}
\innerProd{v}{\Phi'(0)^{-1}(\Phi'(t \xi) - \Phi'(0)) \xi}/t
\right)
\]
whenever $h = \Phi'(0) \xi$. The identity
$g'(0) \Phi'(0) \xi = 2 r \innerProd{v}{\xi}$
% $\| g'(0) \| = 2 r \| \Phi'(0)^{-1} \|$
and \ref{e:prop:StrongConvexityDiffeomorph:L1}
for all $\xi \perp v$, $\xi \not= 0$ imply
\ref{e:lem:C11QuadSupport:b} for all
$h \in \ker g'(0) \setminus \{ 0 \}$, and an application of Lemma
\ref{lem:C11QuadSupport:b} proves our claim.

Now, since the origin is also a boundary point of
$\Phi( \Omega )$, and since the right hand side of
\ref{e:prop:StrongConvexityDiffeomorph:proof:1} is weakly supported at
the origin locally, so is $\Phi( \Omega )$, and hence,
$\interior ( \Phi( \Omega ) )$.
\end{proof}

\subsection{Convexity of attainable sets of discrete-time systems}
\label{ss:convex:DiscreteTime}

We next present a result that enables us to verify hypothesis
\ref{h:quantizer}, and hence, to establish the correctness of the
algorithm proposed in section \ref{s:computation} for the computation
of discrete abstractions of the quantized system
\ref{e:TimeDiscreteAutonomousControlSystem:G},\ref{e:TimeDiscreteAutonomousControlSystem:Q},
whenever a particular quantizer
% \ref{e:TimeDiscreteAutonomousControlSystem:Q}
together with its system $C'$ of operating range cells is given.
In what follows, $D_i^j f$ denotes the partial derivative of
order $j$ with respect to the $i$th argument, of the map $f$.

\begin{theorem}
\label{th:ConvexAttainDiscreteTime}
Assume \ref{h:GdiffeomorphismC1} with smoothness $C^{1,1}$, let
$\psi$ denote the general solution of
\ref{e:TimeDiscreteAutonomousControlSystem:G}, and let
$N \in \mathbb{N}$ and $\Omega \subseteq X$ be $r$-convex with
$\Omega \not= \mathbb{R}^n$. Assume that there are
$L_1,L_2 \in \mathbb{R}$ such that
\begin{align}
\label{e:th:ConvexAttainDiscreteTime:L1}
L_1
& \geq
\alpha_{+}( D_1 G(x,w) )^2 / \alpha_{-}( D_1 G(x,w) ),\\
\label{e:th:ConvexAttainDiscreteTime:L2}
L_2
& \geq
\limsup_{h \to 0}
\frac{\| D_1 G(x,w)^{-1} D_1 G(x + h,w) - \id\|}{\| h\|}
\end{align}
for all
$(x,w) \in \psi(\intco{0;N},\Omega,U^{\mathbb{Z}_{+}}) \times U
\subseteq X \times U$, where $\alpha_{+}(A)$ and
$\alpha_{-}(A)$ denote the maximum and minimum, respectively, singular
values of $A$. Then the attainable set $\psi(k,\Omega,u)$ is convex
for all $k \in \intcc{0;N}$ and all $u \colon \mathbb{Z}_{+} \to U$ if
\begin{equation}
\label{e:th:ConvexAttainDiscreteTime:r}
r L_2 \sum_{\tau = 0}^{N-1} L_1^{\tau}
\leq 1.
\end{equation}
\end{theorem}

\begin{proof}
We may assume $\Omega$ contains at least two points as well as $k = N$
without loss of generality.
By our hypotheses on the right hand side $G$ of
\ref{e:TimeDiscreteAutonomousControlSystem:G}, the map 
$\Phi \defas \psi(N,\cdot,u)$ is a $C^{1,1}$-diffeomorphism between an
open neighborhood of $\Omega$ and an open subset of $X$.
We first prove the claim under the assumption that
\ref{e:th:ConvexAttainDiscreteTime:r} is strict by applying
Prop.~\ref{prop:StrongConvexityDiffeomorph} to $\Phi$:

Let $x \in \boundary \Omega$, $v$ any unit length normal
to $\Omega$ at $x$, and $\xi \perp v$.
% , $\xi \not= 0$.
For $t$ small enough and $k \in \intcc{0;N}$ define
$y_k(t) = D_2 \psi(k, x + t \xi,u) \xi$. Then $y_0(t) = \xi$, and the
sequence $y(t)$ solves the variational equation to
\ref{e:TimeDiscreteAutonomousControlSystem:G} along
$\psi(\cdot, x + t \xi,u)$, i.e.,
\begin{equation}
\label{e:th:ConvexAttainDiscreteTime:proof:ve}
y_{k+1}(t)
=
D_1 G( \psi( k, x + t \xi, u), u_k ) y_k(t)
\end{equation}
for all $k \in \intco{0;N}$.
Next define
$z_k(t) = ( y_k(t) - y_k(0) )/t$ for $t > 0$ small enough. Then
$z_0(t) = 0$, and the sequence $z(t)$ solves another linear difference equation, namely,
\begin{equation}
\label{e:th:ConvexAttainDiscreteTime:proof:vee}
z_{k+1}(t)
=
D_1 G( \psi(k,x,u), u_k ) z_k(t) + b_k
\end{equation}
for all $k \in \intco{0;N}$, where $b_k$ denotes
\[
%\frac{D_1 G( \psi(k, x + t \xi, u), u_k ) - D_1 G( \psi(k, x, u),  u_k )}{t}
%y_k(t).
\left(
D_1 G( \psi(k, x + t \xi, u), u_k ) - D_1 G( \psi(k, x, u),  u_k )
\right)
y_k(t)/t.
\]
Note that in the case of $C^2$-smoothness, if we let $t$ tend to $0$,
then \ref{e:th:ConvexAttainDiscreteTime:proof:vee} reduces to the
variational equation (whose solution is
$D_2^2 \psi(k, x,u) \xi^2$)
of the variational equation
\ref{e:th:ConvexAttainDiscreteTime:proof:ve}.
Now observe
$(k,k_0) \mapsto D_2 \psi(k,x,u) D_2 \psi(k_0,x,u)^{-1}$ is the
transition matrix of the homogeneous system associated with
\ref{e:th:ConvexAttainDiscreteTime:proof:vee}, use the identity
\begin{equation}
\label{e:th:ConvexAttainDiscreteTime:proof:ident}
\Phi'(x)^{-1}
(
\Phi'(x + t \xi)
-
\Phi'(x)
) \xi / t
=
D_2 \psi(N,x,u)^{-1} z_N(t)
\end{equation}
and apply the discrete variation of constants formula \cite{LakshmikanthamTrigiante02}
%solution formula \ref{e:ODifferenceElinGenSol}
to \ref{e:th:ConvexAttainDiscreteTime:proof:vee} to see that the left
hand side of \ref{e:th:ConvexAttainDiscreteTime:proof:ident} equals
\begin{equation}
\label{e:th:ConvexAttainDiscreteTime:proof:ident2}
\sum_{\tau = 0}^{N-1}
D_2 \psi(\tau,x,u)^{-1}
D_1 G( \psi(\tau,x,u), u_\tau )^{-1}
b_\tau.
\end{equation}
From the variational equation of
\ref{e:TimeDiscreteAutonomousControlSystem:G} along
$\psi(\cdot, x, u)$ we obtain
%the estimates
\begin{align*}
\| D_2 \psi( \tau + 1, x, u) \|
&\leq
\alpha_{+}
( D_1 G(p, u_\tau) )
\| D_2 \psi( \tau, x, u) \|,\\
\| D_2 \psi( \tau + 1, x, u)^{-1} \|
&\leq
\frac{\| D_2 \psi( \tau, x, u)^{-1} \|}%
{\alpha_{-}
( D_1 G(p, u_\tau) )
},
\end{align*}
where $p = \psi(\tau,x,u)$; hence,
by our hypothesis \ref{e:th:ConvexAttainDiscreteTime:L1},
\begin{equation}
\label{e:th:ConvexAttainDiscreteTime:proof:est2}
\| D_2 \psi( \tau, x, u)^{-1} \|
\cdot
\| D_2 \psi( \tau, x, u) \|^2
\leq
L_1^{\tau}
\end{equation}
for all $x \in \Omega$, $u \colon \mathbb{Z}_{+} \to U$ and
$\tau \in \intco{0;N}$.

Let $\varepsilon > 0$ be arbitrary.
Then
$\| D_2 \psi(\tau, x + t \xi, u) \| \leq (1 + \varepsilon) \| D_2 \psi(\tau, x, u) \|$ 
whenever $t$ is small enough. Use this fact, the bound
\ref{e:th:ConvexAttainDiscreteTime:L2}, the mean value theorem, and
\ref{e:th:ConvexAttainDiscreteTime:proof:est2} to obtain the upper
bound
$
(1 + \varepsilon)^3 L_2 \| \xi \|^2
\sum_{\tau = 0}^{N-1}
L_1^{\tau}
$
for the norm of \ref{e:th:ConvexAttainDiscreteTime:proof:ident2},
for all $x \in \Omega$, $u \colon \mathbb{Z}_{+} \to U$ and all $t > 0$ small enough.
Now let $\varepsilon$ tend to $0$ to see that the strict variant of
\ref{e:th:ConvexAttainDiscreteTime:r} implies
\ref{e:prop:StrongConvexityDiffeomorph:L1}, hence the convexity of
$\Phi(\Omega)$.

To complete the proof, assume $\Omega$ is of form
\ref{e:def:StrongConvexity} and define
\begin{equation}
\label{e:lem:StrConvexHausdorffConvergence}
\Theta(s)
=
\bigcap_{x \in M}
\cBall(x,s)
\end{equation}
for $s > 0$. Then
$\Phi(\Theta(s))$ is convex for all $s < r$ by the first part of this
proof. By Lemma \ref{lem:StrConvexHausdorffConvergence}, $\Theta(s)$
converges to $\Omega$ in Hausdorff distance, and that property is
preserved under diffeomorphisms. Consequently, $\Phi(\Omega)$ is the
limit of convex sets, and thus, is itself convex \cite{Webster94}.
\end{proof}

We remark that the hypotheses of Theorem \ref{th:ConvexAttainDiscreteTime}
can be verified by inspection of the right hand side $G$ of
\ref{e:TimeDiscreteAutonomousControlSystem:G}. Indeed, suitable
constants $L_1$ and $L_2$ are obtained from estimates of singular
values
of $D_1 G(x,w)$
and of a Lipschitz constant of
$D_1 G(x,w)^{-1} D_1 G(y,w)$ with respect to $y$, respectively.
In this regard, note also that
$L_2$ is just a bound on
$\| D_1 G(x,w)^{-1} D_1^2 G(x,w) \|$ if the right hand side $G$ is of
class $C^2$ with respect to its first argument.

\subsection{Construction of admissible quantizers}
\label{ss:convex:AdmissibleQuantizers}

We now turn to the question of how to construct an admissible
quantizer.
Let there be given some $N \in \mathbb{N}$ and a compact subset
$K \subseteq X$ of the state space $X$ of the discrete-time system
\ref{e:TimeDiscreteAutonomousControlSystem:G} together with an open
neighborhood $V \subseteq X$ of $K$. Intuitively, $N$ is the memory
span of an abstraction we seek to compute, $K$ is the intended
operating range of the quantizer, and $V$ can be thought of as a
maximal operating range, i.e., $X \setminus V$ should be covered by
overflow symbols. Of course, our choice of $N$, $K$ and $V$ would
depend on particularities of the problem we intend to solve.
% ; see also the discussion in Section \ref{s:computation}.

In addition, let there be given a full-dimensional, convex
polyhedron $Y \subseteq \mathbb{R}^n$ together with a strongly
convex set $\widehat Y \not= \mathbb{R}^n$ containing
$Y$ as a subset.
Denote the general solution of the discrete-time system
\ref{e:TimeDiscreteAutonomousControlSystem:G}  by $\psi$.
We will first choose a finite set $C'$ of scaled and translated copies
of $Y$ that cover $K$, i.e., 
\begin{align}
\label{e:eC'}
C'
&=
\Menge{z + \lambda Y}{z \in Z},\\
\label{e:eC'2}
K
&\subseteq
Z + \lambda Y
\end{align}
for some $\lambda > 0$ and some finite subset
$Z \subseteq \mathbb{R}^n$, and then supplement
$C'$ by overflow symbols to obtain a finite cover $C$ of
$\mathbb{R}^n$ for which
\ref{h:quantizer} holds. Our choice will further guarantee
$Z + \lambda \widehat Y \subseteq V$ and that the attainable sets
$\psi(k,z + \lambda \widehat Y, u)$ are convex for all
$k \in \intcc{1;N}$, $z \in Z$, and
$u \colon \intco{0;k} \to U$.

Assume without loss of generality that the origin is an interior point of $Y$ and
$\widehat Y$ is $1/2$-convex, and denote the distance between $K$ and
$\mathbb{R}^n \setminus V$ by $d$.
Then $d > 0$ as $K$ is compact and $V$ is open, and
$K + \lambda \widehat Y$ is compact and contained in $V$ for any
$\lambda \in \intoc{0,d}$. Hence
$K' \defas \psi(\intco{0;N},K + \lambda \widehat Y,U^{\mathbb{Z}_{+}})$ is
compact if $U$ is finite.
If \ref{h:GdiffeomorphismC1} holds with smoothness $C^{1,1}$, then
$G(\cdot,w)$ is a $C^{1,1}$-diffeomorphism, so we may choose
$L_1$, $L_2$ and $\lambda > 0$ such that
$\lambda / 2 \leq r \defas ( L_2 \sum_{\tau = 0}^{N - 1} L_1^{\tau})^{-1}$ and
\ref{e:th:ConvexAttainDiscreteTime:L1} and
\ref{e:th:ConvexAttainDiscreteTime:L2} hold for all
$(x,w) \in K' \times U$, as well as
$K + \lambda \widehat Y \subseteq V$.
Then $\lambda \widehat Y$ is $r$-convex, thus
attainable sets $\psi(k,z + \lambda \widehat Y, u)$ are convex
for all $k \in \intcc{1;N}$, $z \in K$ and
$u \colon \intco{0;k} \to U$ by Theorem
\ref{th:ConvexAttainDiscreteTime}.

Since $\lambda > 0$, $K$ is compact and the origin is an interior
point of $Y$, we can find a finite subset $Z \subseteq K$ for
which \ref{e:eC'2} holds, and we could even guarantee
$K \subseteq Z + \lambda \interior Y$ if necessary.
Finally, Lemma \ref{lem:PolyhedralCovering} in the Appendix shows that
if the set $C'$ is defined by \ref{e:eC'}, it can be supplemented by
convex polyhedra to obtain a finite covering of $\mathbb{R}^n$.
We have thus proved the following result, which easily extends to the
case of a compact rather than finite input alphabet.

\begin{theorem}
\label{cor:th:ConvexAttainDiscreteTime}
Let the input alphabet $U$ of the discrete-time system
\ref{e:TimeDiscreteAutonomousControlSystem:G} be finite and assume
\ref{h:GdiffeomorphismC1} with smoothness $C^{1,1}$.
Let further be given some $N \in \mathbb{N}$, a compact subset
$K \subseteq X$, an open neighborhood $V \subseteq X$ of $K$, as well
as a full-dimensional, convex polyhedron $Y \subseteq \mathbb{R}^n$
together with a strongly convex set
$\widehat Y \subseteq \mathbb{R}^n$ for which
$Y \subseteq \widehat Y \not= \mathbb{R}^n$.\\
Then there is a finite subset $Z \subseteq \mathbb{R}^n$, some
$\lambda > 0$, and a superset $C$ of the set $C'$ defined by
\ref{e:eC'} such that \ref{h:quantizer} and
\[
K
\subseteq
Z + \lambda Y
\subseteq
Z + \lambda \widehat Y
\subseteq
V
\]
hold. Moreover, $\Delta \cap \interior \Delta' = \emptyset$ for all
$\Delta \in C \setminus C'$ and all $\Delta' \in C'$, and one may
additionally require $\Delta \cap K = \emptyset$ for all overflow
symbols $\Delta \in C \setminus C'$.
\end{theorem}

\subsection{Convexity of attainable sets of sampled systems}
\label{ss:convex:Sampled}

We finally provide two results useful for sampled systems.

\begin{theorem}
\label{th:ConvexAttainContinuousTime}
Assume \ref{h:FstetigC1LoesungenAufRplus}, let the right hand side $F$
of \ref{e:TimeContinuousAutonomousControlSystem:F} be of class
$C^{1,1}$ with respect to its first argument, and let
$\varphi$ denote the general solution of
\ref{e:TimeContinuousAutonomousControlSystem:F}. Let
$t > 0$ and $\Omega \subseteq X$ be $r$-convex with
$\Omega \not= \mathbb{R}^n$. Further
assume that there are $M_1,M_2 \in \mathbb{R}$ such that
\begin{align*}
M_1
& \geq
2 \mu_{+}( D_1 F(x,w) ) - \mu_{-}( D_1 F(x,w) ),\\
M_2
& \geq
\limsup_{h \to 0}
\frac{\| D_1 F( x + h, w) - D_1 F(x,w) \|}{\| h\|}
\end{align*}
for all
$(x,w) \in \varphi(\intcc{0,t},\Omega,\mathcal{V}) \times V \subseteq X \times V$,
where $\mu_{+}(A)$ and $\mu_{-}(A)$ denote the maximum and minimum,
respectively, eigenvalues of the symmetric part $(A + A^{\ast})/2$ of
$A$. Then the attainable set $\varphi(\tau,\Omega,v)$ is convex for
all $\tau \in \intcc{0,t}$ and all admissible input signals
$v \in \mathcal{V}$ if
\begin{equation}
\label{e:th:ConvexAttainContinuousTime:r}
r M_2 \int_{0}^{t} \exp \left( M_1 \rho \right) d \rho
\leq 1.
\end{equation}
\end{theorem}

\begin{proof}
We may assume $\Omega$ contains at least two points as well as $\tau = t$
without loss of generality.
By our hypotheses on the right hand side $F$ of
\ref{e:TimeContinuousAutonomousControlSystem:F}, the map 
$\Phi \defas \varphi(t,\cdot,v)$ is a $C^{1,1}$-diffeomorphism between
an open neighborhood of $\Omega$ and an open subset of $X$.
We assume $\Omega$ is of form
\ref{e:def:StrongConvexity}, define $\Theta(s)$ for $s > 0$ by
\ref{e:lem:StrConvexHausdorffConvergence}, and
prove $\varphi(t,\Theta(s),v)$ is convex for any
$s \in \intoo{0,r}$ by applying
Proposition \ref{prop:StrongConvexityDiffeomorph} to $\Phi$.
The theorem then follows from Lemma
\ref{lem:StrConvexHausdorffConvergence}.

To this end, set $I = \intcc{0,t}$ and
$X' = \varphi(I,\interior \Omega, v)$, and define
$f \colon I \times X' \to \mathbb{R}^n$ by
$f(\tau,x) = F(x,v(\tau))$. Since $X'$ is an open neighborhood of
$\varphi(I,\Theta(s),v)$, 
% ODE \ref{e:ODifferentialE} with right hand
% side $f$
the ODE $\dot x = f(t, x)$
fulfills the hypothesis of \cite[Theorem 3]{i07Convex}. The
proof of the latter result shows that if $v$ is continuous, then for
any $x \in X'$ and any $\varepsilon > 0$ we have
$\| \Phi'(x)^{-1} ( \Phi'(x + h) - \Phi'(x) ) \|
\leq (1 + \varepsilon)^2 M_2 \| h \| \int_{0}^{t} \exp \left( M_1 \rho \right) d \rho$
for all sufficiently small $h$.
The extension to piecewise continuous $v$ is straightforward.
Then \ref{e:th:ConvexAttainContinuousTime:r} implies
\ref{e:prop:StrongConvexityDiffeomorph:L1} with $s$ substituted for
$r$ and $\zeta$ substituted for $v$,
for any $x \in \boundary \Theta(s)$ and any $\zeta,\xi \in \mathbb{R}^n$
with $\| \zeta \| = 1$ and $\xi \not= 0$. Thus $\Phi(\Theta(s))$ is convex
by Proposition \ref{prop:StrongConvexityDiffeomorph}.
\end{proof}

\begin{theorem}
\label{th:C2useful}
Assume \ref{h:FstetigC1LoesungenAufRplus}, let the right hand side $F$
of \ref{e:TimeContinuousAutonomousControlSystem:F} be of class
$C^2$ with respect to its first argument,
and let $\varphi$, $t$, $\Omega$, and $r$ be as in Theorem
\ref{th:ConvexAttainContinuousTime}. Assume further that there is a constant
$L_2$ such that
{\multlinegap0pt
\begin{multline}
\label{e:th:C2useful:L2}
\left\| \int_0^{\delta}
D_2 \varphi(\tau,x,v)^{-1}
D_1^2 F( p(\tau), v( \tau ) )
\left( D_2 \varphi(\tau,x,v) h \right)^2
d\tau
\right\|
\ifCLASSOPTIONonecolumn\relax\else\\\fi
\leq
L_2
\| h \|^2
\end{multline}}%
for all
$x \in \Omega$,
$\delta \in \intcc{0,t}$,
$v \in \mathcal{V}$ and $h \in \mathbb{R}^n$, where
$p(\tau) = \varphi(\tau,x,v)$.
Then the attainable set $\varphi(\tau,\Omega,v)$ is convex for all
$\tau \in \intcc{0,t}$ and all admissible input signals
$v \in \mathcal{V}$ if $r L_2 \leq 1$.
\end{theorem}

\begin{proof}
As in the proof of Theorem \ref{th:ConvexAttainContinuousTime}, we
assume $\Omega$ contains at least two points and $\tau = t$, define
$\Theta(s)$ by \ref{e:lem:StrConvexHausdorffConvergence}, and observe
$\Phi \defas \varphi(t,\cdot,v)$ is a $C^2$-diffeomorphism.
Let $x \in X$, $h \in \mathbb{R}^n$, $v \in \mathcal{V}$, and
define $y(t) = D_2 \varphi(t,x,v) h$. Then $y(0) = h$, and $y$
solves the variational equation to
\ref{e:TimeContinuousAutonomousControlSystem:F} along
$\varphi(\cdot,x,v)$, i.e.,
\begin{equation}
\label{e:th:C2useful:ve}
\dot y(t) =
D_1 F (\varphi(t,x,v), v(t)) y(t)
\end{equation}
for all $t \geq 0$. Next define
$z(t) = D_2^2 \varphi(t,x,v) h^2$. Then $z(0) = 0$, and $z$
solves another linear ODE, namely,
\begin{equation}
\label{e:th:C2useful:vee}
\dot z(t) =
D_1 F (\varphi(t,x,v), v(t)) z(t)
+
D_1^2 F (\varphi(t,x,v), v(t)) y(t)^2
\end{equation}
for all $t \geq 0$. (Note that $x$ is a parameter rather than an
initial value in \ref{e:th:C2useful:ve}.)
Now observe
$(t,t_0) \mapsto D_2 \varphi(t,x,v) D_2 \varphi(t_0,x,v)^{-1}$ is the
transition matrix of the homogeneous system associated with
\ref{e:th:C2useful:vee} and apply the solution
formula for linear differential equations \cite{Hartman02}
to \ref{e:th:C2useful:vee} to see that
$\Phi'(x)^{-1} \Phi''(x) h^2$ equals the integral in
\ref{e:th:C2useful:L2} with $t$ substituted for $\delta$.
Proposition \ref{prop:StrongConvexityDiffeomorph} shows
$\Phi(\Theta(s))$ is convex, and the theorem follows from Lemma
\ref{lem:StrConvexHausdorffConvergence}.
\end{proof}

Theorems \ref{th:ConvexAttainContinuousTime} and \ref{th:C2useful} with the choice
$t = N\cdot T$
provide sufficient conditions for attainable sets of the sampled
system \ref{e:TimeDiscreteAutonomousControlSystem:G} to be convex, as
required in hypothesis \ref{h:quantizer} of Section \ref{s:computation}.
In the case of Theorem \ref{th:ConvexAttainContinuousTime}, that condition
can be verified directly from
properties of the right hand side $F$ of the continuous-time system
\ref{e:TimeContinuousAutonomousControlSystem:F} by estimating
eigenvalues and a Lipschitz constant, with $M_2$ being just a
bound on $\| D_1^2 F(x,w) \|$ if $F$ is of class $C^2$ with respect to
its first argument.
In contrast, application of Theorem \ref{th:C2useful} requires
estimating the integrand in \ref{e:th:C2useful:L2}. That higher effort
often pays off when it results in larger bounds on $r$ than
\ref{e:th:ConvexAttainContinuousTime:r}. Note that, in view of the algorithm
proposed in this paper, larger bounds will typically translate
into lower computational complexity; see
\ref{fig:GeometricIntuitionConvexityForSmallTimes}(b).

\section{Example}
\label{s:ex}
\begin{figure*}[!t]
\centering
% \psfrag{0}[r][r]{$x_1 = 0$}
% \psfrag{p}[r][r]{$x_1 = \pi$}
\psfrag{0}[r][r]{}
\psfrag{pi}[r][r]{}
\psfrag{x1}[][]{$x_1$}
\psfrag{u}[][]{$u$}
\psfrag{p}[r][r]{}
\psfrag{-->}[r][r]{}
\psfrag{Mg}[b][b]{}
\psfrag{Lg}[br][br]{}
\psfrag{Rg}[bl][bl]{}
\includegraphics[width=.495\linewidth]{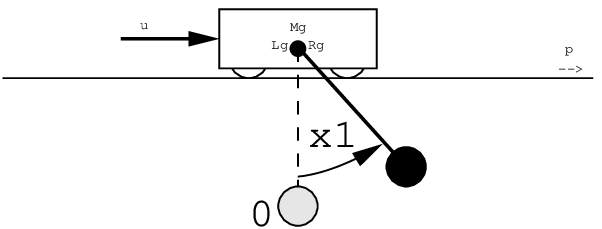}
\psfrag{pla}[][]{pendulum}
\psfrag{u1}[][]{$u_1$}
\psfrag{superv}[][]{supervisor}
\psfrag{discr}[][]{discrete}
\psfrag{act}[][]{actuator}
\psfrag{gen}[][]{\ifCLASSOPTIONonecolumn\small\fi{}generator}
\psfrag{low}[][]{low-level}
\psfrag{con}[][]{controller}
\psfrag{u}[][]{$u$}
\psfrag{x}[][]{$x$}
\psfrag{Om}[lb][lb]{}
\includegraphics[width=.495\linewidth]{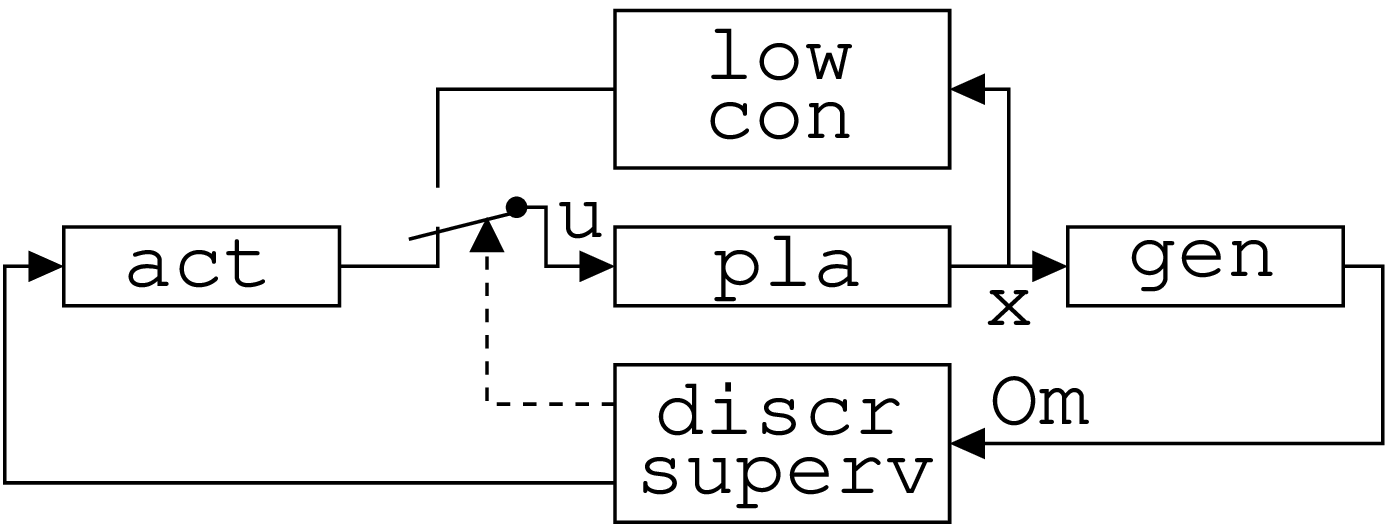}\\
\hspace*{\fill}(a)\hspace*{\fill}\hspace*{\fill}(b)\hspace*{\fill}
\caption{\label{fig:PendulumAndHybridControlLoop}
Example investigated in Section \ref{s:ex}:
The pendulum (a) is swung up by the hybrid control system (b).}
\end{figure*}

In this section we shall demonstrate an application of our results
from Sections \ref{s:computation} and \ref{s:convex} within the
framework of abstraction based supervisory control of sampled systems
by solving a nonlinear, global problem with constraints.
To begin with, consider the system
\begin{subequations}
\label{e:Pendulum}
\begin{align}
\dot x_1 &= x_2,\\
\dot x_2 &= -\omega^2 \sin(x_1) - u\; \omega^2 \cos(x_1) - 2 \gamma x_2,
\end{align}
\end{subequations}
which describes the motion of a pendulum mounted on a cart.
Here, $\omega$ and $\gamma$ are parameters, specifically, $\gamma$ is
a friction coefficient, and $x_1$ is the angle
between the pendulum and the downward vertical. See
\ref{fig:PendulumAndHybridControlLoop}(a).
\begin{figure*}[t]
\centering
%% hier beginnt (a) %%
\psfrag{1}[b][b]{\small$\geq 1$}
\psfrag{2}[b][b]{\small$\geq 2$}
\psfrag{3}[b][b]{\small$\vphantom{\geq} 3$}
\psfrag{x0}[l][b]{\small$\,0$}
\psfrag{x1}[r][b]{\small$1\,$}
\psfrag{x2}[bl][b]{\small$2$}
\psfrag{x3}[bl][b]{\small$3$}
\psfrag{x4}[r][b]{\small$4\,$}
\psfrag{x5}[bl][b]{}
\psfrag{x6}[b][b]{\small$6$}
\psfrag{x7}[l][b]{}
\psfrag{x8}[r][b]{\small$8\,$}
\psfrag{x9}[b][b]{}
\psfrag{x10}[b][b]{\small $10$}
\psfrag{x11}[br][b]{\small $11$}
\psfrag{x12}[b][b]{\small $12$}
\psfrag{x13}[b][b]{}
\psfrag{x14}[bl][b]{\small $14$}
\psfrag{x15}[b][b]{}
\psfrag{x16}[rt][]{\small $16\,$}
\psfrag{x17}[b][b]{}
\psfrag{x18}[b][b]{\small $18$}
\psfrag{x19}[l][b]{}%\small $19$}
\psfrag{x20}[bl][b]{\small\color{white}$20$}
\psfrag{0}[][]{\small$0$}
\psfrag{p1}[][]{\small$\vphantom{0}\pi$}
\psfrag{p2}[r][]{\small$2\pi$}
\psfrag{yp}[r][r]{\small$\pi$}
\psfrag{ym}[t][r]{\tiny$-$\small$\pi\;\,$}
\psfrag{E}[][]{\Large\color{blue}$\mathbf{E}$}
\psfrag{S}[r][]{\large S}
\psfrag{Z}[][]{\Huge\color{lightgray}Z}
\psfrag{MemorySpan}[t][t]{\small memory span}
\includegraphics[width=.495\linewidth]{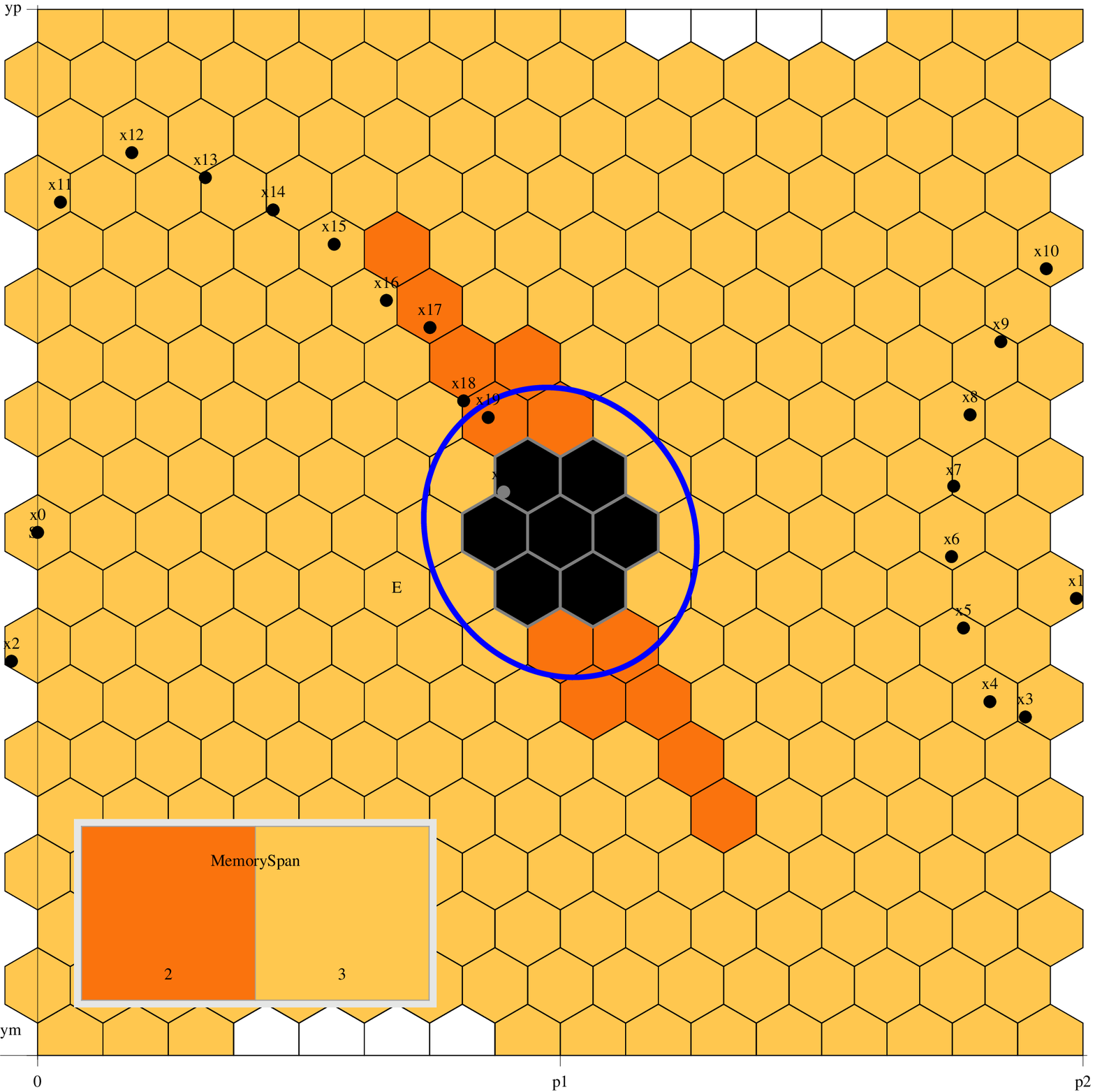}
%% hier beginnt (b) %%
\psfrag{1}[b][b]{\small$\geq 1$}
\psfrag{2}[b][b]{\small$\geq 2$}
\psfrag{3}[b][b]{\small$\vphantom{\geq} 3$}
\psfrag{x0}[l][b]{\small$\,0$}
\psfrag{x1}[lt][b]{\small$\,1$}
\psfrag{x2}[][]{}
\psfrag{x3}[bl][b]{\small$3$}
\psfrag{x4}[][]{}
\psfrag{x5}[bl][b]{\small$5$}
\psfrag{x6}[][]{}
\psfrag{x7}[lt][b]{\small$\,7$}
\psfrag{x8}[][]{}
\psfrag{x9}[b][b]{\small$9$}
\psfrag{x10}[b][b]{}
\psfrag{x11}[b][b]{\small$\,11$}
\psfrag{x12}[b][b]{\small$12$}
\psfrag{x13}[b][b]{}
\psfrag{x14}[b][b]{\small$14$}
\psfrag{x15}[b][b]{}
\psfrag{x16}[rt][]{}
\psfrag{x17}[b][b]{}
\psfrag{x18}[tl][b]{\small$18$}
\psfrag{x19}[b][b]{}
\psfrag{x20}[b][b]{\small $20$}
\psfrag{x21}[b][b]{}
\psfrag{x22}[b][b]{}
\psfrag{x23}[b][b]{}
\psfrag{x24}[b][b]{\small$24\;$}
\psfrag{x25}[b][b]{}
\psfrag{x26}[b][b]{}
\psfrag{x27}[b][b]{}
\psfrag{x28}[b][b]{\small$28$}
\psfrag{x29}[b][b]{\small$29\,$}
\psfrag{x30}[b][b]{}
\psfrag{x31}[b][b]{}
\psfrag{x32}[b][b]{\small$\;32$}
\psfrag{x33}[b][b]{}
\psfrag{x34}[b][b]{}
\psfrag{x35}[tr][b]{\small$35\,$}
\psfrag{x36}[b][b]{}
\psfrag{x37}[b][b]{}
\psfrag{x38}[t][b]{\raisebox{-1.1em}{\small\color{white}$38$}}
\psfrag{0}[][]{\small$0$}
\psfrag{p1}[][]{\small$\vphantom{0}\pi$}
\psfrag{p2}[r][]{\small$2\pi$}
\psfrag{yp}[r][r]{\small$\pi$}
\psfrag{ym}[t][r]{\tiny$-$\small$\pi\;\,$}
%\psfrag{E}[tr][tr]{\color{blue}$\mathbf{E}$}
\psfrag{E}[][]{\Large\color{blue}$\mathbf{E}$}
\psfrag{S}[r][]{\large S}
\psfrag{Z}[][]{\Huge\color{lightgray}Z}
\psfrag{H}[][]{\Huge\color{lightgray}H}
\psfrag{MemorySpan}[t][t]{\small memory span}
\includegraphics[width=.495\linewidth]{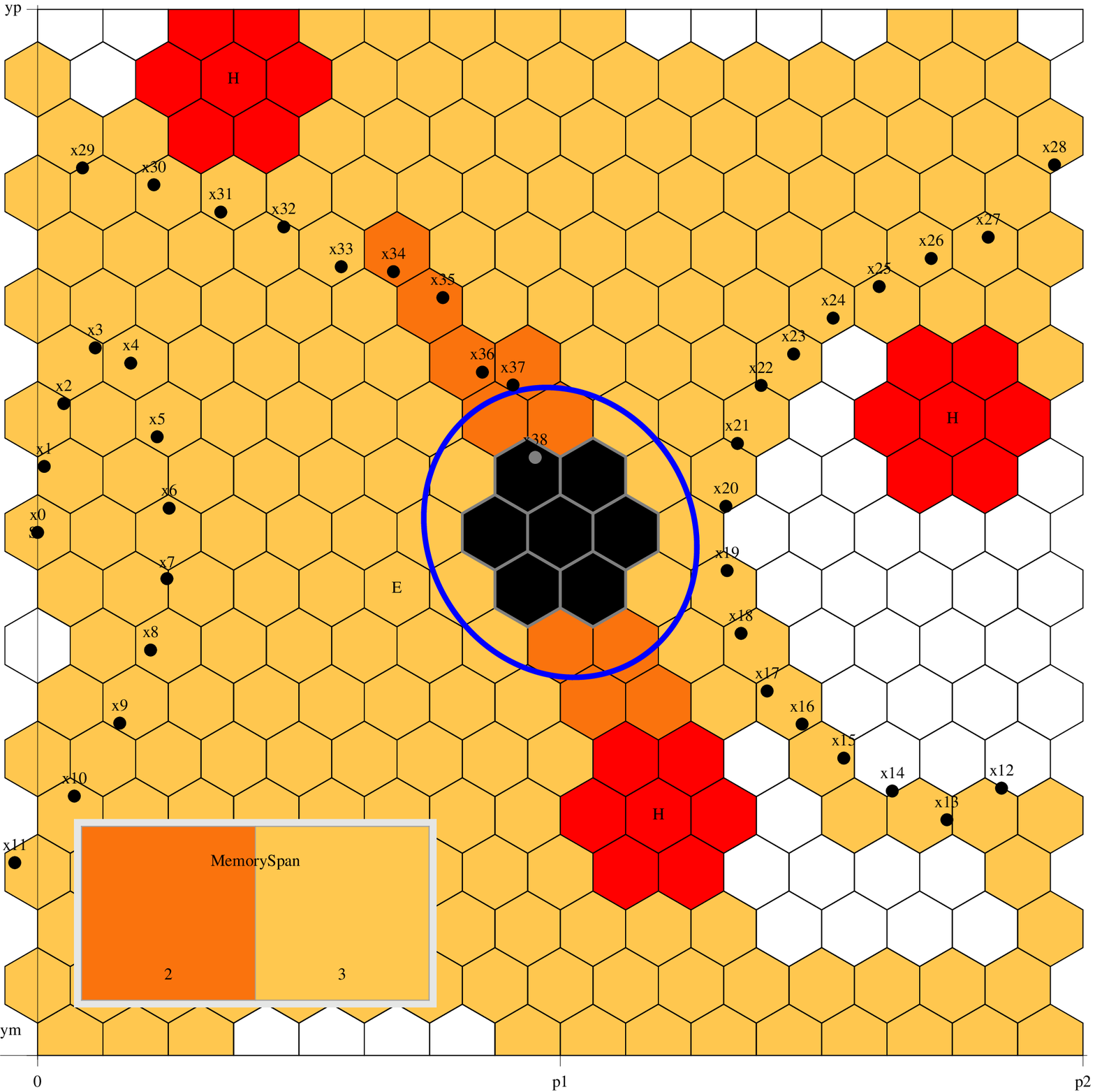}\\
\hspace*{\fill}(a)\hspace*{\fill}\hspace*{\fill}(b)\hspace*{\fill}
\caption{\label{fig:PendulumSupervisorResults}
(a) The state space of the pendulum system \ref{e:Pendulum} is covered
with quantizer cells.
Supervisors designed on the basis of abstractions of memory span
$N \in \{2,3\}$ force the sampled and quantized pendulum system into
the stability region $E$ of the low-level controller, from anywhere in
the indicated regions. The region for $N = 3$ contains the origin;
one particular trajectory is shown.
(b) Problem from (a) with extra constraints in the form of three
obstacles in state space, which are labeled $H$ in the illustration.
}
\end{figure*}%
The motion of the cart is not modeled; its acceleration $u$ is
considered
a control.

We seek to swing up the pendulum by means of the hybrid control system
shown in \ref{fig:PendulumAndHybridControlLoop}(b), which possesses a simple
hierarchical structure.
The low-level controller
%$u_1 \colon \mathbb{R}^2 \to \mathbb{R}$
%$u_1 \colon x \mapsto u$
is to stabilize the pendulum at its upright position.
That is, the point $(\pi,0)$ becomes an asymptotically stable equilibrium
of the closed loop composed of \ref{e:Pendulum} and the low-level
controller, so there will be some non-trivial, positively invariant
subset $E$ of its stability region.
The supervisor, on the other hand,
%is a discrete, possibly dynamic controller that
would force the state
% of \ref{e:Pendulum}
from some neighborhood of
the origin into $E$
% within finite time
and on success, would hand over control to the low-level controller.
The supervisor will be realized by a finite
automaton, which is why it is connected to the continuous plant via
interface devices \cite{KoutsoukosAntsaklisStiverLemmon00},
to the effect that the open loop composed of actuator, pendulum and
generator is represented by
the sampled and quantized system
\ref{e:TimeDiscreteAutonomousControlSystem:G},\ref{e:TimeDiscreteAutonomousControlSystem:Q}
associated with \ref{e:Pendulum}.

A suitable low-level controller together with a positively invariant
set $E$ is straightforward to determine
\cite{Sontag98,HenrionGarulli05}. For example,
if $\omega > 0$ and $0 \leq \gamma \leq \omega$,
the affine state feedback
$u = 2 ( \pi - x_1 - x_2 / \omega )$ stabilizes \ref{e:Pendulum} at $(\pi,0)$,
with the positively invariant ellipsoid
\[
E
=
(\pi,0)
+
\Menge{x \in \mathbb{R}^2}%
  { 63 \omega^2 x_1^2 + 12 \omega x_2 x_1 + 56 x_2^2 \leq 42 \omega^2 }
\]
being a subset of the stability region.
Here we focus on the design of quantizer and supervisor.
So, let us consider a sampled version
\ref{e:TimeDiscreteAutonomousControlSystem:G} of \ref{e:Pendulum}
and impose the constraints
\begin{subequations}
\label{e:StateAndControlConstraints}
\begin{align}
\label{e:StateConstraint}
| x_2 | & \leq \pi,\\
\label{e:ControlConstraint}
| u | & \leq 2,
\end{align}
\end{subequations}
which model physical limits of an experimental setup.
For the sake of simplicity, we choose controls to be constant on a
common sampling interval $\intcc{0,T}$, specifically,
\[
U = \{
t \mapsto 0,
t \mapsto -2,
t \mapsto 2
\}
\]
in the notation of Section \ref{ss:SystemClass}. This guarantees the
sampled system \ref{e:TimeDiscreteAutonomousControlSystem:G} fulfills
hypothesis \ref{h:GdiffeomorphismC1}.

We next design a suitable quantizer
\ref{e:TimeDiscreteAutonomousControlSystem:Q} as part of the generator
device in the hybrid control system of
\ref{fig:PendulumAndHybridControlLoop}(b). To this end, assume that
the control constraint \ref{e:ControlConstraint} is satisfied and that
problem data and sampling period are given by
\[
\omega = 1,
\gamma = 0.01,
T = 0.2.
\]
Theorem \ref{th:pendel} in the Appendix shows that the attainable
set $\varphi(t,\Omega,u)$ is convex whenever
$\Omega \subseteq \mathbb{R}^2$ is $r$-convex, $r > 0.4$, and
$0 \leq t \leq 3 T = 0.6$, where $\varphi$ denotes the general
solution of \ref{e:Pendulum}.
This implies any translated and possibly truncated copy of the regular
hexagon given by its set
\begin{equation}
\label{e:unithexagon}
\frac{\pi}{16 \sqrt{3}}
\{
(0, \pm 2), (\sqrt{3}, \pm 1), (-\sqrt{3}, \pm 1)
\}
\end{equation}
of vertices,
which has circumradius $\pi / (8 \sqrt{3}) < 0.23$, may be chosen as a
quantizer cell in the computation of abstractions of memory span up to
$3$. Further, in view of the state constraint \ref{e:StateConstraint} we
may restrict our investigation of the dynamics of \ref{e:Pendulum} to the
region $K$ defined by
$
K= \mathbb{R} \times \intcc{-\pi,\pi}
$.
So, let us choose $C'$ as a set of $304$ translated copies of the
hexagon \ref{e:unithexagon}, each intersected with $K$. This intersection
either leaves a hexagon unchanged or results in an irregular
pentagon; see \ref{fig:PendulumSupervisorResults}.
Finally, supplement $C'$ with two overflow symbols,
\[
C = C' \cup \{ \mathbb{R} \times \intco{\pi,\infty},\mathbb{R} \times
\intoc{-\infty,-\pi} \}.
\]

Note that since the right hand side of \ref{e:Pendulum} is periodic in $x$ with
period $(2 \pi,0)$, we have implicitly considered the system \ref{e:Pendulum} on
the cylinder
\cite{Hartman02,Sontag98}.
Having said this, $C$ can really be regarded as a covering of the
state space of \ref{e:TimeDiscreteAutonomousControlSystem:G}.

With the choices we have made above, hypotheses \ref{h:GdiffeomorphismC1} and
\ref{h:quantizer} in Sections \ref{s:prelims} and
\ref{s:computation} are fulfilled.
In particular, Theorem \ref{th:pendel} in the Appendix shows that for
each cell $\Delta \in C'$, we may choose
the smallest intersection of six closed balls of radius $0.4$
containing $\Delta$ for the set $\widehat \Delta$ in hypothesis
\ref{h:quantizer};
see \ref{fig:HexagonAndStronglyConvexHullAndTheirImage}.
We finally choose $\Sigma(\widehat \Delta)$
to consist of six supporting half-spaces of $\widehat \Delta$
as in \ref{fig:HexagonAndStronglyConvexHullAndTheirImage},
and analogously for the pentagons in $C'$.
The set $\Sigma(\widehat \Delta)$ is then supplied to the algorithm in
\ref{alg:DiscreteAbstractions} to compute abstractions of the sampled
and quantized pendulum system defined earlier. The results are
summarized in \ref{tab:Statistics}:
\begin{table}[t]
\caption{\label{tab:Statistics}
Computation of abstractions of the pendulum example.% in Section \ref{s:ex}.%
}
\centering
\begin{tabular}{crrrr}
N & half-spaces & polyhedra & states & transitions\\
\hline
1 & 7170  & 41059  & 306   & 4246  \\
2 & 22914 & 97203  & 4552  & 35734 \\
3 & 69048 & 351523 & 36040 & 220442
\end{tabular}
\end{table}
The memory span $N$ of the abstractions we have computed,
the number of half-spaces determined from solutions of ODE
\ref{e:TimeContinuousAutonomousControlSystem:F:adjoint}, the number of
polyhedra tested for emptiness, and the number of states and transitions
in a finite automaton realization \cite{Willems89,Moor99,MoorRaisch99}
of the abstraction. 
The data in \ref{tab:Statistics} highlights the fact that
half-spaces are shared among polyhedra, which is an important feature
of the algorithm from Section \ref{s:computation}.
Indeed, while each transition corresponds to a non-empty
intersection of half-spaces, the number of those transitions by far
exceeds the total number of computed half-spaces if $N > 1$.

In order to obtain a suitable supervisor for the control system of
\ref{fig:PendulumAndHybridControlLoop}(b), we solve certain auxiliary
control problems posed in terms of the abstractions already computed.
So, let
$N \in \{1,2,3\}$, denote by $B_N$ the abstraction of memory span $N$
that we have computed, define a start region $S$ and a target region $Z$ by
\begin{align*}
S
&=
\Menge{ \Delta \in C' }{ (0,0) \in \Delta },\\
Z
&=
\Menge{ \Delta \in C' }{ \Delta \subseteq E },
\end{align*}
see \ref{fig:PendulumSupervisorResults},
and consider the following problem:
Determine the supervisor in the form of a map
%\begin{equation}
%\label{e:Supervisor}
$R
\colon
\bigcup_{k = 1}^{N}
U^{k - 1}
\times
C^{k}
\to
U
$
%\end{equation}
such that whenever $(u,\Delta) \in B_N$,
$\Delta_0 \in S$, and
$
%\[
%\begin{equation}
%\label{e:BehaviorOfSupervisor}
%\forall_{k \in \mathbb{Z}_{+}}
u_k
=
R
(
u|_{\intoo{k - N; k}},
\Delta|_{\intoc{k - N; k}}
)
%\end{equation}
%\]
$
for all $k \in \mathbb{Z}_{+}$,
then there is some $k \in \mathbb{Z}_{+}$ such that
$
%\[
\Delta_k \in Z
$ and $
%\forall_{\tau \in \intco{0;k}}
\Delta_\tau \in C'
%\]
$ for all $\tau \in \intco{0;k}$.
This specification requires that if $(u,\Delta)$ is any
signal that may possibly be produced by the closed loop composed of
the supervisor $R$ and some plant that realizes the behavior $B_N$,
and if that signal additionally starts in $S$, then it remains in $C'$
until it eventually enters $Z$.
While its specification is not a complete behavior \cite{Willems89},
that kind of discrete control problem is equivalent to a shortest path
problem in some hypergraph
\cite{OezverenWillskyAntsaklis91},
and thus, can be efficiently solved
\cite{GalloLongoPallotinoNguyen93}.
In fact, we have been able to obtain a supervisor $R$ in the case
$N = 3$, and to prove there is none enforcing the above
specification if $N \in \{ 1, 2 \}$, with run times negligible
compared to the ones observed in the computation of the abstractions.
Specifically, for $N = 3$, the target region is reached within at
most $27$ steps.
It follows that $R$ is compatible with the actual
plant, i.e., with the sampled and quantized pendulum system defined
earlier, and also enforces the above
specification when combined with that plant rather than
with the abstraction $B_3$, on which the design of $R$ was based
\cite{MoorDavorenAnderson02}. See also
\ref{fig:PendulumSupervisorResults}(a).

\begin{figure*}[t]
%\centering
%\hspace*{\fill}
%% fuer Graphik am linken Rand:
%
\psfrag{y}[b][b]{cpu}
\psfrag{y00}[r][r]{\footnotesize$1$}
\psfrag{y10}[r][r]{\footnotesize$10$}
\psfrag{y20}[r][r]{\footnotesize$10^2$}
\psfrag{y30}[r][r]{\footnotesize$10^3$}
%% hier beginnt (a) %%
\psfrag{x}[l][l]{\footnotesize$|U|$}
\psfrag{2}[t][t]{\footnotesize$2$}
\psfrag{3}[t][t]{\footnotesize$3$}
\psfrag{5}[t][t]{\footnotesize$5$}
\psfrag{7}[t][t]{\footnotesize$7$}
\psfrag{9}[t][t]{\footnotesize$9$}
\psfrag{11}[t][t]{\footnotesize$11$}
\psfrag{21}[t][t]{\footnotesize$21$}
\psfrag{41}[t][t]{\footnotesize$41$}
\includegraphics[width=.32775\linewidth]{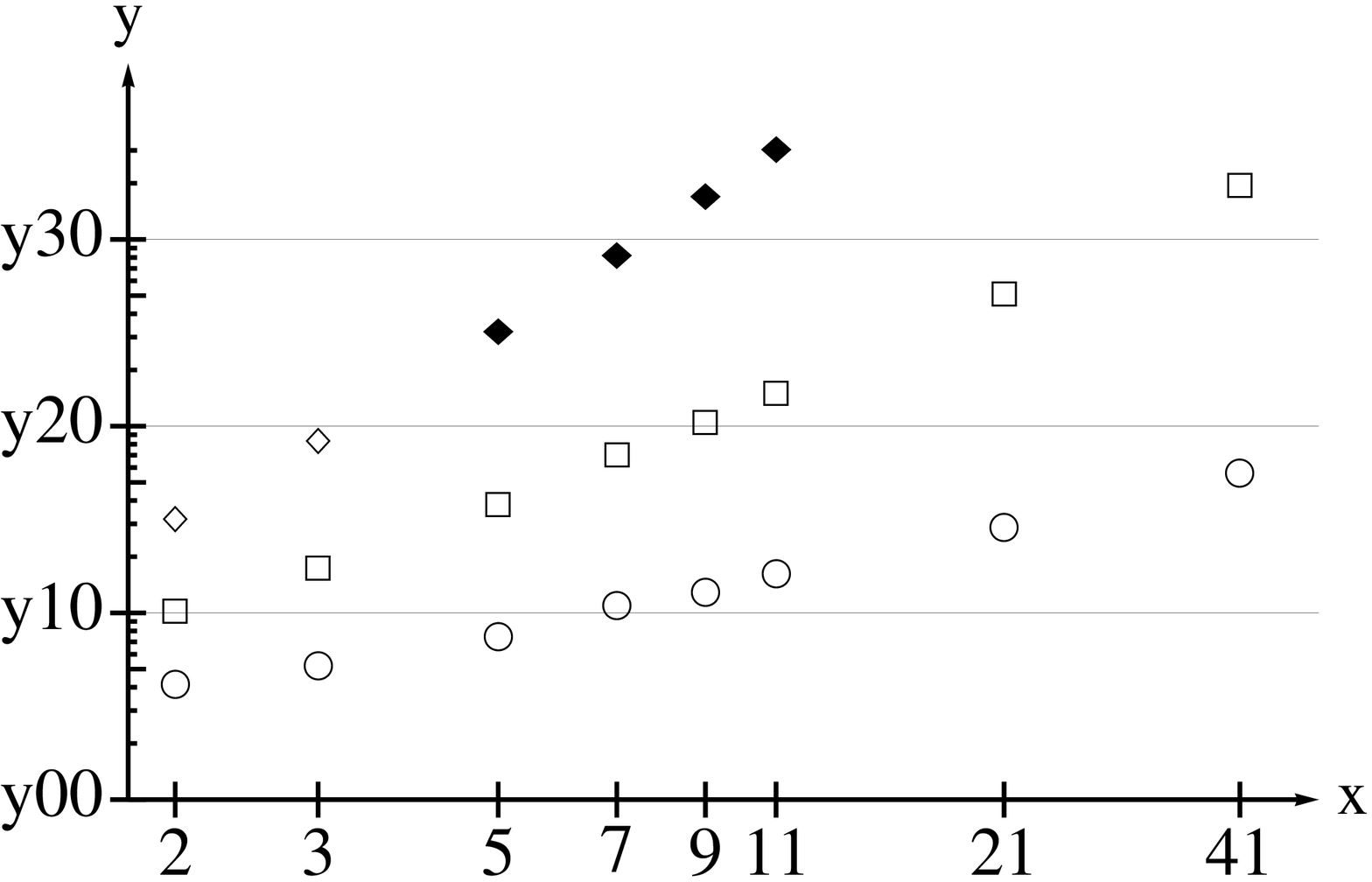}%
%\hspace*{\fill}
% fuer alle weiteren Graphiken
\psfrag{y00}[r][r]{}
\psfrag{y10}[r][r]{}
\psfrag{y20}[r][r]{}
\psfrag{y30}[r][r]{}
%% hier beginnt (b) %%
\psfrag{x}[l][l]{\footnotesize$m$}
\psfrag{6}[t][t]{\footnotesize$6$}
\psfrag{12}[t][t]{\footnotesize$12$}
\psfrag{24}[t][t]{\footnotesize$24$}
\psfrag{48}[t][t]{\footnotesize$48$}
\psfrag{90}[t][t]{\footnotesize$90$}
\includegraphics[width=.32775\linewidth]{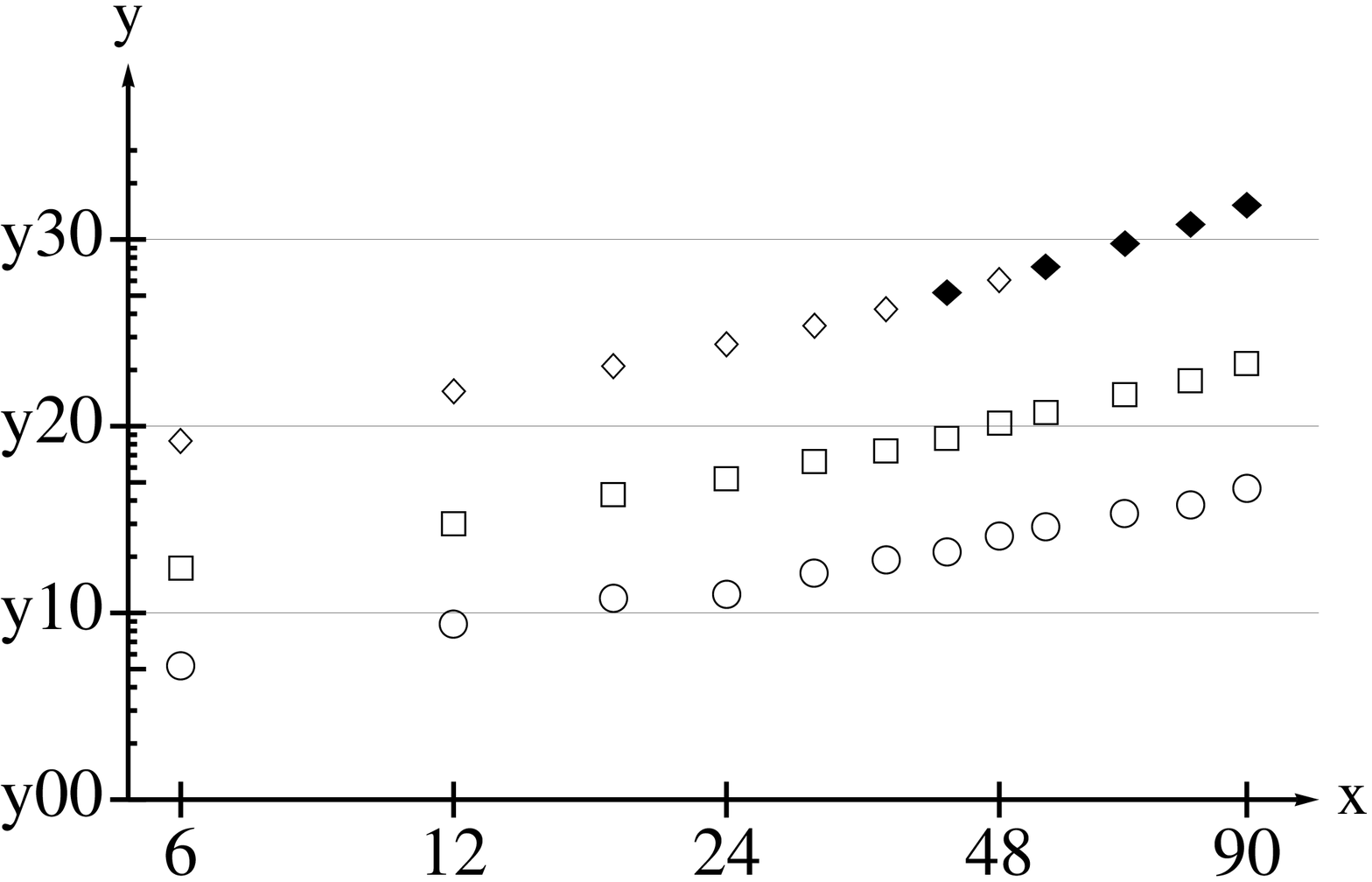}%
%\hspace*{\fill}
%% hier beginnt (c) %%
\psfrag{x}[l][l]{\footnotesize$|C|$}
\psfrag{200}[t][t]{\footnotesize$200$}
\psfrag{500}[t][t]{\footnotesize$500$}
\psfrag{x30}[t][t]{\footnotesize$1000$}
\psfrag{x40}[t][t]{\footnotesize$10000$}
\includegraphics[width=.32775\linewidth]{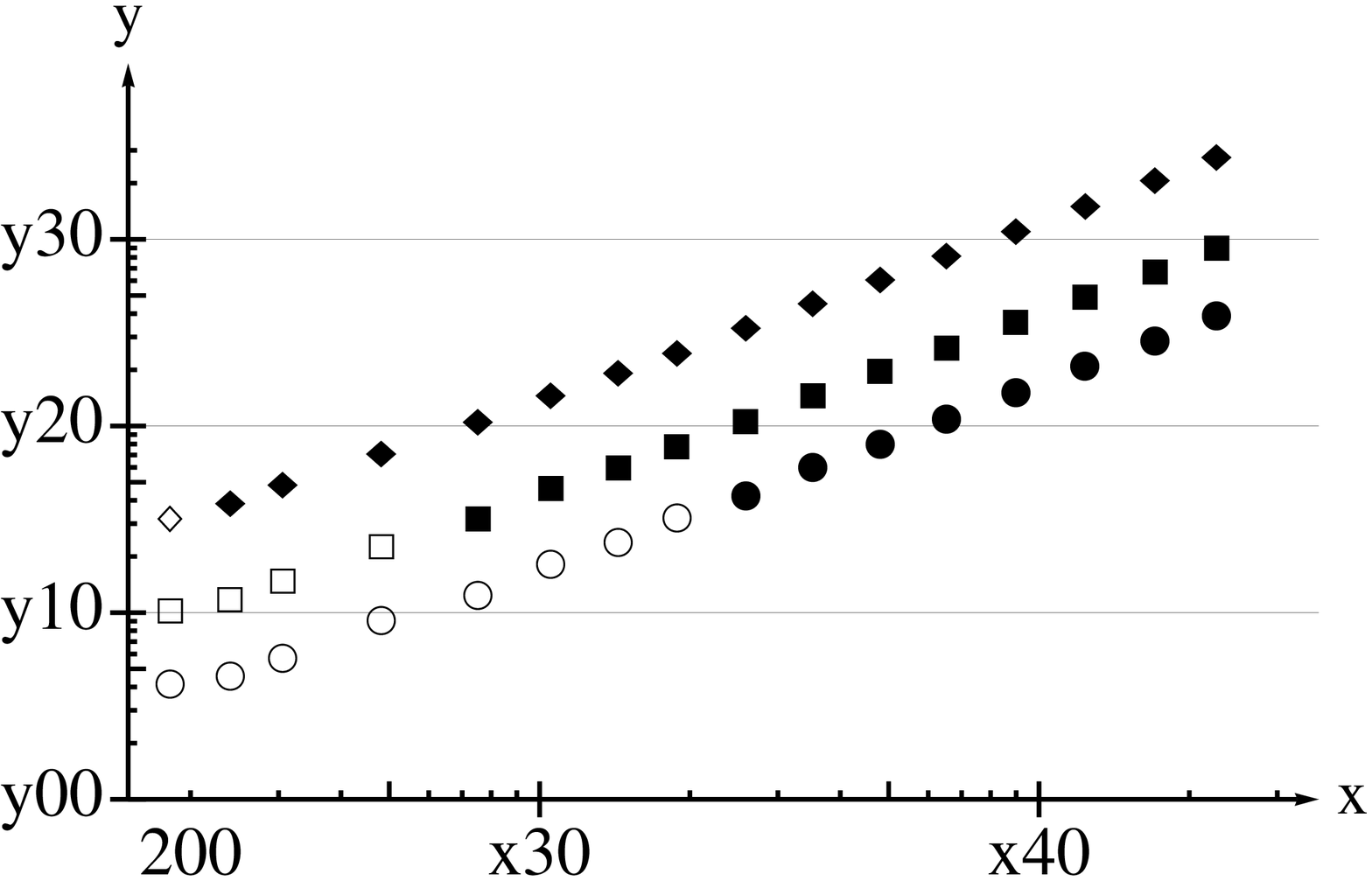}
\hspace*{\fill}\\
\hspace*{\fill}(a)\hspace*{\fill}\hspace*{\fill}(b)\hspace*{\fill}\hspace*{\fill}(c)\hspace*{\fill}
\caption{\label{fig:ComputationalResults}
Dependence of run times in seconds on the number $|U|$ of controls, the
number $m$ of half-spaces supporting $\widehat \Delta$, and the number
$|C|$ of quantizer cells, of an implementation of the algorithm in
\ref{alg:DiscreteAbstractions} with \begriff{Mathematica 7.0}
\cite{Mma5}, run on four threads of an Intel Xeon CPU E$5620$
($2.4$ GHz). Results for memory spans $1$, $2$ and $3$ correspond to
the symbols $\circ$, {\tiny$\square$}, and {\scriptsize$\lozenge$},
respectively, which are filled iff the corresponding
abstraction has lead to a solution of the control problem considered in
Section \ref{s:ex}.
(a) $|C| = 182$, $m = 6$.
(b) $|U| = 3$, $|C| = 182$.
(c) $|U| = 2$, $m = 6$.
}
\end{figure*}%
Modifying the above example,
we have varied the number of controls, the
number of supporting hyperplanes, and the number of quantizer
cells. See \ref{fig:ComputationalResults}.
Given the fact that the reported run times have been obtained from
interpreted rather than from compiled code, we expect that the
algorithm in \ref{alg:DiscreteAbstractions} can also be successfully
applied to systems essentially more complex than \ref{e:Pendulum}. See
also \cite{i10product}. \ref{fig:ComputationalResults}(b) also
demonstrates the importance of accurately approximating attainable
sets. In particular, we have verified for the quantizer used in
\ref{fig:ComputationalResults}(b) that the control problem considered
in this section could not be solved using ellipsoidal rather than
strongly convex supersets of quantizer cells.

We have additionally
investigated a scenario with extra constraints in the form of
obstacles in state space, by simply treating the obstacle cells as
overflow symbols. The results illustrated in
\ref{fig:PendulumSupervisorResults}(b) show that our approach
would also be feasible in the presence of complicated constraints such
as the ones regularly met in motion planning problems
\cite{LaValle06}.

Finally, we would like to point out that the supervisors we have designed solve
control problems for a sampled version of \ref{e:Pendulum}.
In fact, the discrete abstractions we are using are not capable of
representing the evolution of continuous-time systems between sampling
times, and hence, control problems for the latter systems cannot be
treated directly.
However, solutions for an important class of continuous-time control
problems can be obtained from solutions of auxiliary problems for
sampled systems using a robust version of the original specification,
e.g. \cite{FainekosGirardKressGazitPappas09}.
For the problem considered in this section, a robust
specification could easily be obtained by tightening state
constraints, i.e., by decreasing the bound \ref{e:StateConstraint} and
enlarging the obstacles in \ref{fig:PendulumSupervisorResults}(b) by a
suitable amount.

\section{Conclusions}
\label{s:Conclusions}

We have presented a novel algorithm for the computation of discrete
abstractions of nonlinear systems as well as a set of sufficient
conditions for the convexity of attainable sets. While the usefulness
of the first relies on the second contribution, the latter may be of
separate interest. Practicability of our results in the design of
discrete controllers for nonlinear continuous plants under state and
control constraints has been demonstrated by an example, and we also
expect their use to be of advantage in attainability and verification
problems.

The algorithm proposed in this paper is the first one that not
only yields abstractions of finite but otherwise arbitrary memory span
that are suitable for solving general control problems, but also applies to
nonlinear
% discrete-time and sampled
systems under rather mild
conditions, which essentially reduce to sufficient smoothness in the
case of sampled systems. Previous approaches are
confined to abstractions of memory span $1$ with two exceptions, which
apply only to monotone dynamics \cite{MoorRaisch02} and are not
rigorous and limited to
solving reachability problems \cite{GrueneMueller08}, respectively.
We emphasize that increasing the memory span
may be the only way to improve the accuracy of abstractions up to a level at
which analysis and synthesis problems can be solved. One example is
networked control systems \cite{MatveevSavkin09}, where quantization
effects are part of the systems to be investigated and state space
quantizations cannot be arbitrarily refined.

Its wide applicability and the fact that it builds on relatively
simple computations distinguishes our approach from competing
techniques even if we restrict ourselves to abstractions of memory
span $1$. In particular, some methods only apply to systems whose
continuous-valued dynamics is defined by ordinary difference and
differential equations with multi-affine or polynomial right hand
sides \cite{MalerBatt08,BermanHalaszKumar07,GirardPappas05}, or
require stability \cite{Tabuada08,PolaTabuada09} or
deriving a state space partition in accordance with the exact
system dynamics prior to their application
\cite{Broucke98,CainesWei98,StiverKoutsoukosAntsaklis01}.
Others require the use of interval arithmetic
\cite{JaulinWalter97,StursbergKowalewskiEngell00,AlthoffStursbergBuss09,TazakiImura09},
deciding satisfiability of formulas over certain logical theories \cite{Tiwari08b},
or solving complex optimization problems
\cite{ChutinanKrogh03,MitchellBayenTomlin05,GrueneJunge07}.

Finally, a distinctive feature of our method is that the error by
which attainable sets of quantizer cells are over-approximated is
quadratic in the size of the cells.
This has been achieved by extending our earlier results from
\cite{i07Convex,i07MMAR} to apply to strongly convex sets rather than
merely ellipsoids.
Due to this improvement, our method will outperform the sampling method
\cite{Junge99} and any other method whose approximation error depends
linearly on the dispersion \cite{LaValle06} of some grid, e.g.
\cite{Tabuada08,PolaTabuada09,i09HSCC}, whenever highly accurate abstractions
are to be computed.

The techniques we propose can currently be applied to systems with
finite input alphabets only and additionally depend on the
ability to design suitable quantizers. The latter can be quite
demanding, despite the results presented here. An extension to systems
with continuous inputs and an automated procedure for designing
quantizers would considerably enhance our method.
It should also be extended to account for disturbances and
uncertainties, including numerical discretization errors and the
effects of finite arithmetic in order to address robustness issues
and to obtain validated results.

\section*{Acknowledgment}
\noindent
The author thanks B.~Farrell (M\"unchen) and
O.~Stursberg (Kassel) for their valuable comments on earlier versions
of this paper, as well as A.~Herrmann (Berlin), for providing
convenient access to computing facilities.

\appendix
\label{s:app}

\begin{lemma}
\label{lem:sConvexUnitBall}
Let $r > 0$, $z \in \mathbb{R}^n$,
$x,y \in \cBall(z,r)$, $x \not= y$, and
$s = \| x - y \|^2/(8r)$.
Then $\cBall((x+y)/2,s) \subseteq \cBall(z,r)$.
\end{lemma}

\begin{proof}
Let $A$ be an arc of a circle of radius $r$ joining $x$ and $y$ whose
length does not exceed $\pi r$. Then $A \subseteq \cBall(z,r)$,
e.g. \cite{FrankowskaOlech81}, and
$\min \Menge{ \| a - (x + y)/2 \|}{a \in A} = r - (r^2 - \| x - y \|^2/4)^{1/2}$.
The latter is easily shown to be bounded below by $s$.
\end{proof}

The \begriff{Hausdorff distance} between any non-empty, compact
subsets $M, N \subseteq \mathbb{R}^n$ is defined to be the infimum of
$r \in \mathbb{R}_{+}$ for which both $M \subseteq N + \cBall(0,r)$
and $N \subseteq M + \cBall(0,r)$ \cite{Webster94}.

\begin{lemma}
\label{lem:StrConvexHausdorffConvergence}
Let $M \subseteq \mathbb{R}^n$, $M \not= \emptyset$,
and define
$
\Theta(s)
=
\bigcap_{x \in M}
\cBall(x,s)
$
for all $s > 0$.
If $r > 0$ and $\Theta(r)$ contains at least two points, then
$
\lim_{s \to r, s < r} \Theta(s) = \Theta(r)
$
in Hausdorff distance.
\end{lemma}

\begin{proof}
$\Theta(s)$ is convex and compact for any $s > 0$, and
$\Theta(r)$ possesses nonempty interior by Lemma
\ref{lem:sConvexUnitBall}; hence
$\Theta(r) = \closure ( \interior ( \Theta(r) ) )$.
If $p \in \interior \Theta(r)$, then $p \in \Theta(s)$ whenever $s$ is
sufficiently close to $r$, in particular,
$\Theta(s) \not= \emptyset$. Thus the Hausdorff distance between
$\Theta(s)$ and $\Theta(r)$ is well-defined. Moreover, given
$\varepsilon > 0$ and $y \in \Theta(r)$, there is
$p \in \interior \Theta(r)$ with $\| p - y \| < \varepsilon$, hence
$y \in \Theta(s) + \oBall(0,\varepsilon)$ for some $s < r$, where
$\oBall(x,r)$ denotes the open Euclidean ball of radius $r$ centered
at $x$. This shows
$\Theta(r) \subseteq \bigcup_{s \in \intoo{0,r}} (\Theta(s) + \oBall(0,\varepsilon))$,
and compactness of $\Theta(r)$ implies
$\Theta(r) \subseteq \Theta(s) + \oBall(0,\varepsilon)$ for some $s < r$,
hence for all $s < r$ sufficiently close to $r$.
\end{proof}

\begin{lemma}
\label{lem:PolyhedralCovering}
Let $P_1, \dots, P_k \subseteq \mathbb{R}^n$ be convex polyhedra. Then
$\closure ( \mathbb{R}^n \setminus \bigcup_{i = 1}^k P_i )$ is the
finite union of convex polyhedra.
\end{lemma}

\begin{proof}
For any polyhedron
$M = \Menge{ x \in \mathbb{R}^n}{ A x \leq b}$,
$A$ an $m \times n$-matrix, $b \in \mathbb{R}^m$, $m \geq 1$, the closure of
$\mathbb{R}^n \setminus M$ equals
$\bigcup_{j = 1}^m \Menge{ x \in \mathbb{R}^n}{ A_j x \geq b_j}$,
where $A_j$ denotes the $j$th row of $A$. Therefore,
%\[
$
\closure
\left(
\mathbb{R}^n \setminus \bigcup_{i = 1}^k P_i
\right)
=
\bigcap_{i = 1}^k \closure \left( \mathbb{R}^n \setminus P_i \right)
=
\bigcap_{i = 1}^k
\bigcup_{j = 1}^m
Q_{i,j}
$
%\]
for suitable half-spaces $Q_{i,j}$.
Since intersection and union distribute over each other, the right hand side
of the previous identity equals
$\bigcup_{j \in J^I}
\bigcap_{i \in I}
Q_{i,j_i}
$,
where $I = \intcc{1;k}$ and $J = \intcc{1;m}$.
\end{proof}

\begin{theorem}
\label{th:pendel}
Let $t > 0$ and assume the input $u$ to the pendulum equations
\ref{e:Pendulum} is piecewise continuous with
$|u(\tau)| \leq \hat u$ for all $\tau \in [0,t]$.
Define
\begin{align*}
\widehat \omega
&=
\max\left\{1,|\omega| \left( 1 + \hat u^2 \right)^{1/4}\right\},\\
r&=
\frac{12 \widehat \omega^2 \;(1+(\widehat\omega  +\gamma)^2)^{-3/2}}%
{
 \sinh(3 \widehat \omega t) +
 \sinh(\widehat \omega t)
 \left(
%   \frac{12 \hat \omega^3}{(\hat \omega^2 + 1)^{3/2}}-3
   12 (\widehat \omega^{-2} + 1)^{-3/2} - 3
 \right)},
\end{align*}
where $\max$ denotes the maximum,
and assume $0 \leq \gamma \leq \frac{3}{4} \widehat \omega$ and
%$2\kappa_{-} T \leq \pi$.
$2 ( \widehat \omega^2 - \gamma^2 )^{1/2} t \leq \pi$.
Then the attainable set $\varphi(t,\Omega,u)$ is convex for any
$r$-convex subset $\Omega \subseteq \mathbb{R}^2$,
where $\varphi$ denotes the general solution of \ref{e:Pendulum}.
\end{theorem}%

\begin{proof}
First note that the right hand side of \ref{e:Pendulum} is linearly
bounded \cite{Hartman02}, which implies
$\varphi(\tau,\mathbb{R}^2,u) = \mathbb{R}^2$
for any $\tau \in \mathbb{R}_{+}$. One may therefore assume
$\Omega \not= \mathbb{R}^2$ without loss of generality.
Then apply Theorem \ref{th:C2useful} and use
the estimate
\[
\omega^2
| u(\tau) \cos(\varphi(\tau,x_0,u)_1) + \sin(\varphi(\tau,x_0,u)_1) |
\leq \widehat \omega^2,
\]
the fact that $D_2 \varphi(\cdot,x,u)$ fulfills the variational
equation to \ref{e:Pendulum} along $\varphi(\cdot,x,u)$, Cramer's
rule, and the formula of Abel--Liouville \cite{Hartman02}
to see that it suffices to show
\begin{equation}
\label{e:th:SufficientConvexityConditionPendulum:EuclideanBall:b:2}
\widehat \omega^2
\int_0^{t}
  \e^{2\gamma\tau}
  \left\|( D_2\varphi(\tau,x_0,u))_{1,\cdot} \right\|^3
d\tau
\leq
1/r.
\end{equation}
Here, the subscript ``$1,\cdot$'' denotes the first row.
Now set
$\kappa = \left(\widehat \omega^2 + \gamma^2 \right)^{1/2}$
and observe
\ifCLASSOPTIONonecolumn\linebreak\fi
$
\left( 1 + (\kappa + \gamma)^2 \right)
\kappa^{-2}
\leq
\left( 1 + (\widehat \omega + \gamma)^2 \right)
\widehat \omega^{-2}
$
to obtain the upper bound
\begin{equation}
\label{e:th:SufficientConvexityConditionPendulum:EuclideanBall:b:ExpBound:pre}
\frac{\e^{-2\gamma \tau}}{2 \widehat \omega^2}
(1+(\widehat \omega + \gamma)^2)
\left(
\cosh(2 \kappa \tau) + \frac{\widehat \omega^2 - 1}{\widehat \omega^2 + 1}
\right)
\end{equation}
for the squared norm of the first row of $\exp\left(\tau
\left(\begin{smallmatrix}0 & 1\\ \hat \omega^2 & -2 \gamma \end{smallmatrix}\right)
\right)$.
The second step of the proof of \cite[Theorem 6]{i07Convex} shows
\ref{e:th:SufficientConvexityConditionPendulum:EuclideanBall:b:ExpBound:pre}
is also an upper bound for $\| (D_2 \varphi(\tau,x_0,u))_{1,\cdot} \|^2$.
Next show
\begin{equation}
\label{e:th:SufficientConvexityConditionPendulum:EuclideanBall:b:h}
h(\alpha)
\defas
\frac{\cosh\left(( 9 \beta^2 + (\alpha + 4 \beta)^2 )^{1/2}\right)}%
{\cosh(\alpha + 4 \beta)}
-\e^\beta
\leq 0
\end{equation}
for all $\alpha, \beta \geq 0$. The choice
$\alpha = 2(\widehat \omega - 4 \gamma/3)\tau$,
$\beta = 2 \gamma \tau/3$ then yields the upper bound
\begin{equation}
\label{e:th:SufficientConvexityConditionPendulum:EuclideanBall:b:ExpBound}
\frac{\e^{-4\gamma \tau/3}}{\widehat \omega^2}
(1+(\widehat \omega + \gamma)^2)
\left(
\cosh(\widehat \omega \tau)^2 - \frac{1}{\widehat \omega^2 + 1}
\right)
\end{equation}
for $\| (D_2 \varphi(\tau,x_0,u))_{1,\cdot} \|^2$.
Indeed, the map $h$ defined in
\ref{e:th:SufficientConvexityConditionPendulum:EuclideanBall:b:h} is
continuous, and for every $\alpha > 0$, $h'(\alpha)$ exists and is a
positive multiple of
\begin{equation}
\label{e:th:SufficientConvexityConditionPendulum:EuclideanBall:b:TANH}
%\frac{\tanh(\mu + \nu)}{\tanh(\mu)}
\tanh(\mu + \nu)/\tanh(\mu)
-
%\frac{\mu + \nu}{\mu},
(\mu + \nu)/\mu,
\end{equation}
where $\mu = \alpha + 4 \beta$ and
$\nu = (9 \beta^2 + \mu^2)^{1/2} - \mu$.
\ref{e:th:SufficientConvexityConditionPendulum:EuclideanBall:b:TANH}
is monotonically decreasing with respect to $\nu \geq 0$ and vanishes
for $\nu = 0$. Thus $h$ is
monotonically decreasing on $\mathbb{R}_+$. This proves
\ref{e:th:SufficientConvexityConditionPendulum:EuclideanBall:b:h}
since $h(0) \leq 0$.

Finally,
consider the map $g$ defined on $[0,1]$ by
\[
g(s) =
1 - s
\left(
1 - \frac{\widehat \omega^3}{(\widehat \omega^2 + 1)^{3/2}} 
\right)
-
\left(
1-\frac{s}{\widehat \omega^2 + 1}
\right)^{3/2}.
\]
This map is concave since $g''(s) < 0$, and $g(0)=g(1)=0$; thus $g$
is non-negative. For the choice
$s=\cosh(\widehat \omega \tau)^{-2}$, this together with the bound
\ref{e:th:SufficientConvexityConditionPendulum:EuclideanBall:b:ExpBound}
implies
$\| D_2 \varphi(\tau,x_0,u)_{1,\cdot} \|^3$ does not exceed
\def\myequation{%
%\frac{\e^{-2 \gamma \tau}}{\widehat \omega^3}
\e^{-2 \gamma \tau} \widehat \omega^{-3}
(1+(\widehat \omega+\gamma)^2)^{3/2}
\ifCLASSOPTIONonecolumn\relax\else\cdot\\\fi
\left(
\cosh(\widehat \omega \tau)^3
-
\cosh(\widehat \omega \tau)
\left(
%1-\frac{\widehat \omega^3}{(\widehat \omega^2 + 1)^{3/2}}
1 - \widehat \omega^3 (\widehat \omega^2 + 1)^{-3/2}
\right)
\right),
}
\ifCLASSOPTIONonecolumn
\[
\myequation
\]
\else
\begin{multline*}
\myequation
\end{multline*}
\fi%
which directly implies
\ref{e:th:SufficientConvexityConditionPendulum:EuclideanBall:b:2}.
\end{proof}

\bibliographystyle{IEEEtran}
\bibliography{preambles,mrabbrev,strings,fremde,eigeneCONF,eigeneJOURNALS,eigenePATENT,eigeneREPORTS,eigeneTALKS,eigeneTHESES}
\end{document}